\documentclass[11pt]{amsart}%
\usepackage{amsmath}
\usepackage{amsfonts}
\usepackage{amssymb}
\usepackage{graphicx}
\theoremstyle{plain}
\newtheorem{thm}{Theorem}
\newtheorem{cor}[thm]{Corollary}
\newtheorem{lem}[thm]{Lemma}
\newtheorem{prop}[thm]{Proposition}

\newtheorem{example}[thm]{Example}

\begin{document}
\title{Minimality Properties of Tsirelson Type Spaces}

\begin{abstract}
In this paper, we study minimality properties of partly modified mixed
Tsirelson spaces. A Banach space with a normalized basis $(e_{k})$ is said to
be \emph{subsequentially minimal} if for every normalized block basis
$(x_{k})$ of $(e_{k}) ,$ there is a further block $(y_{k})$ of $(x_{k}) $ such
that $(y_{k})$ is equivalent to a subsequence of $(e_{k}).$ Sufficient
conditions are given for a partly modified mixed Tsirelson space to be
subsequentially minimal and connections with Bourgain's $\ell^{1}$-index are
established. It is also shown that a large class of mixed Tsirelson spaces
fails to be subsequentially minimal in a strong sense.

\end{abstract}
\author[D. Kutzarova]{Denka Kutzarova}
\address{Institute of Mathematics, Bulgarian Academy of Sciences and Department of
Mathematics, University of Illinois at Urbana-Champaign, Urbana, IL 61801}
\email{denka@math.uiuc.edu}
\author[D. H. Leung]{Denny H. Leung}
\address{Department of Mathematics, National University of Singapore, 2 Science Drive
2, Singapore 117543}
\email{matlhh@nus.edu.sg}
\author[A. Manoussakis]{Antonis Manoussakis}
\address{Department of Mathematics, University of Aegean, Karlovasi, Samos, GR 83200, Greece}
\email{amanouss@aegean.gr}
\author[W.-K. Tang]{Wee-Kee Tang}
\address{Mathematics and Mathematics Education, National Institute of Education \\
Nanyang Technological University, 1 Nanyang Walk, Singapore 637616}
\email{weekee.tang@nie.edu.sg}
\thanks{The first and second authors are pleased to acknowledge the support they
received to attend the Workshop in Linear Analysis at Texas A \& M University
in August 2005, where the research reported in this paper was begun.}
\maketitle

The class of mixed Tsirelson spaces plays an important role in the structure
theory of Banach spaces and has been well investigated (e.g., \cite{AD, ADKM,
AF, LT, AM,OTW}). In this paper, we will study aspects of the subspace
structure of mixed Tsirelson spaces and (partly) modified mixed Tsirelson
spaces (see definitions below). We are particularly interested in properties
connected with minimality. A infinite-dimensional Banach space $X$ is
\emph{minimal} if every infinite-dimensional subspace has a further subspace
isomorphic to $X$. The work of Gowers \cite{Go} had motivated some recent
studies on minimality (e.g., \cite{F}, \cite{FR}, \cite{P}).

A Banach space $X$ with a normalized basis $\left(  e_{k}\right)  $ is said to
be \emph{subsequentially minimal} if for every normalized block basis $\left(
x_{k}\right)  $ of $\left(  e_{k}\right)  ,$ there is a further block $\left(
y_{k}\right)  $ of $\left(  x_{k}\right)  $ such that $\left(  y_{k}\right)  $
is equivalent to a subsequence of $\left(  e_{k}\right)  .$ It is well known
that the Tsirelson space $T\left[  \left(  {\mathcal{S}}_{1},1/2\right)
\right]  $ has the property that every normalized block basis of its standard
basis is equivalent to a subsequence of $\left(  e_{k}\right)  $ \cite{CJT}.
In particular, it is subsequentially minimal. In \cite[Theorem 9]{LT1}, it was
shown that if a nonincreasing null sequence $\left(  \theta_{n}\right)  $ in
$(0,1)$ is regular ($\theta_{m+n}\geq\theta_{m}\theta_{n}$) and satisfies
\[
(\dag)\quad\quad\lim_{m}\limsup_{n}\frac{\theta_{m+n}}{\theta_{n}}>0,
\]
then the space $T\left[  \left(  {\mathcal{S}}_{n},\theta_{n}\right)
_{n=1}^{\infty}\right]  $ is subsequentially minimal if and only if every
block subspace of $T\left[  \left(  {\mathcal{S}}_{n},\theta_{n}\right)
_{n=1}^{\infty}\right]  $ admits an $\ell^{1}$-$\mathcal{S}_{\omega}%
$-spreading model, if and only if every block subspace of $T\left[  \left(
{\mathcal{S}}_{n},\theta_{n}\right)  _{n=1}^{\infty}\right]  $ has Bourgain
$\ell^{1}$-index greater than $\omega^{\omega}.$ In particular, if $\sup
_{n}\theta_{n}^{1/n}=1,$ then the mixed Tsirelson space $T\left[  \left(
{\mathcal{S}}_{n},\theta_{n}\right)  _{n=1}^{\infty}\right]  $ is
subsequentially minimal \cite{AM}.

This paper is divided into two parts. In the first part, we investigate the
analogs of the results quoted above in the context of partly modified mixed
Tsirelson spaces. In this connection, it is worthwhile to point out that a
subsequentially minimal partly modified mixed Tsirelson space is
quasi-minimal in the sense of Gowers \cite{Go}. Since these spaces are
strongly asymptotic $\ell^{1}$, by \cite{DFKO} they do not contain minimal
subspaces and therefore, they are strictly quasi-minimal. The only typical
known example of a strictly quasi-minimal space was the Tsirelson space. While
that space satisfies the so called blocking principle \cite{CJT}, among our
examples of strictly quasi-minimal spaces there are cases which do not satisfy
that principle. The subsequentially minimal mixed Tsirelson spaces, mentioned
above, are also quasi-minimal, however it is not known if they are strictly
quasi-minimal (see the remarks in \cite{DFKO}). In the second part of the
paper, we give a general sufficient condition for a (unmodified) mixed
Tsirelson space to fail to be subsequentially minimal in a strong sense.

\section{Preliminaries}

Denote by ${\mathbb{N}}$ the set of natural numbers. For any infinite subset
$M$ of ${\mathbb{N}}$, let $[M]$, respectively $[M]^{<\infty}$, be the set of
all infinite and finite subsets of $M$ respectively. These are subspaces of
the power set of ${\mathbb{N}}$, which is identified with $2^{{\mathbb{N}}}$
and endowed with the topology of pointwise convergence. A subset
${\mathcal{F}}$ of $[{\mathbb{N}}]^{<\infty}$ is said to be \emph{hereditary}
if $G\in{\mathcal{F}}$ whenever $G\subseteq F$ and $F\in{\mathcal{F}}$. It is
\emph{spreading} if for all strictly increasing sequences $(m_{i})_{i=1}^{k}$
and $(n_{i})_{i=1}^{k}$, $(n_{i})_{i=1}^{k}\in{\mathcal{F}}$ if $(m_{i}%
)_{i=1}^{k}\in{\mathcal{F}}$ and $m_{i}\leq n_{i}$ for all $i$. We also call
$(n_{i})_{i=1}^{k}$ a spreading of $(m_{i})_{i=1}^{k}$. A \emph{regular}
family is a subset of $[{\mathbb{N}}]^{<\infty}$ that is hereditary, spreading
and compact (as a subspace of $2^{{\mathbb{N}}}$). If $I$ and $J$ are nonempty
finite subsets of ${\mathbb{N}}$, we write $I<J$ to mean $\max I<\min J$. We
also allow that $\emptyset<I$ and $I<\emptyset$. For a singleton $\{n\}$,
$\{n\}<J$ is abbreviated to $n<J$.
If $\mathcal{F},\mathcal{G}\subseteq\lbrack{\mathbb{N}}]^{<\infty}$, let
\[
{\mathcal{F}}[{\mathcal{G}}]=\{\cup_{i=1}^{k}G_{i}:G_{i}\in{\mathcal{G}}%
,G_{1}<\cdots<G_{k},(\min G_{i})_{i=1}^{k}\in\mathcal{F}\}
\]
and
\[
({\mathcal{F}},{\mathcal{G}})=\{F\cup G:F<G,F\in{\mathcal{F}},G\in
{\mathcal{G}}\}.
\]
Inductively, set $({\mathcal{F}})^{1}={\mathcal{F}}$ and $({\mathcal{F}%
})^{n+1}=({\mathcal{F}},({\mathcal{F}})^{n})$ for all $n\in{\mathbb{N}}$. It
is clear that ${\mathcal{F}}[{\mathcal{G}}]$ and $({\mathcal{F}},{\mathcal{G}%
})$ are regular if both ${\mathcal{F}}$ and ${\mathcal{G}}$ are. A class of
regular families that has played a central role is the class of generalized
Schreier families \cite{AA}.

Let ${\mathcal{S}}_{0}$ consist of all singleton subsets of ${\mathbb{N}}$
together with the empty set. Then define ${\mathcal{S}}_{1}$ to be the
collection of all $A\in\lbrack{\mathbb{N}}]^{<\infty}$ such that $|A|\leq\min
A$ together with the empty set, where $|A|$ denotes the cardinality of the set
$A$. If ${\mathcal{S}}_{\alpha}$ has been defined for some countable ordinal
$\alpha$, set ${\mathcal{S}}_{\alpha+1}={\mathcal{S}}_{1}[{\mathcal{S}%
}_{\alpha}]$. For a countable limit ordinal $\alpha$, specify a sequence
$(\alpha_{n})$ that strictly increases to $\alpha$. Then define
\[
{\mathcal{S}}_{\alpha}=\{F:F\in{\mathcal{S}}_{\alpha_{n}}\text{ for some
$n\leq\min F\}\cup\{\emptyset\}$}.
\]
Given a nonempty compact family ${\mathcal{F}}\subseteq\lbrack{\mathbb{N}%
}]^{<\infty}$, let ${\mathcal{F}}^{(0)}={\mathcal{F}}$ and ${\mathcal{F}%
}^{(1)}$ be the set of all limit points of ${\mathcal{F}}$. Continue
inductively to derive ${\mathcal{F}}^{(\alpha+1)}=({\mathcal{F}}^{(\alpha
)})^{(1)}$ for all ordinals $\alpha$ and ${\mathcal{F}}^{(\alpha)}=\cap
_{\beta<\alpha}{\mathcal{F}}^{(\beta)}$ for all limit ordinals $\alpha$. The
\emph{index} $\iota({\mathcal{F}})$ is taken to be the smallest $\alpha$ such
that ${\mathcal{F}}^{(\alpha+1)}=\emptyset$. Since $[{\mathbb{N}}]^{<\infty}$
is countable, $\iota({\mathcal{F}})<\omega_{1}$ for any compact family
${\mathcal{F}}\subseteq\lbrack{\mathbb{N}}]^{<\infty}$. It is well known that
$\iota({\mathcal{S}}_{\alpha})=\omega^{\alpha}$ for all $\alpha<\omega_{1}$
\cite[Proposition 4.10]{AA}.

A sequence $(x_{n})$ in a normed space said to \emph{dominate} a sequence
$(y_{n})$ in a possibly different space if there is a finite constant $K$ such
that $\Vert\sum a_{n}y_{n}\Vert\leq K\Vert\sum a_{n}x_{n}\Vert$ for all
$(a_{n})\in c_{00}$. If two sequences dominate each other, then they are
\emph{equivalent}\textbf{ }and we write $\left(  x_{n}\right)  \sim\left(
y_{n}\right)  $. If $\left(  e_{n}\right)  $ is a basic sequence and
$F\subseteq\mathbb{N}$, $[(e_{n})_{n\in F}]$ denotes the closed linear space
of $\left\{  e_{n}:n\in F\right\}  .$ If $(e_{n})$ is a normalized basis of
$X$, then by $\left(  x_{n}\right)  \prec\left(  e_{n}\right)  $ or $(x_{n})
\prec X$ we shall mean that $\left(  x_{n}\right)  $ is a normalized block
basis of $\left(  e_{n}\right)  .$ We say that $Y$ is a \emph{block subspace}
of $X$, $Y\prec X$, if $X$ has a basis $\left(  x_{n}\right)  $ and $Y=
[(y_{n})_{n\in\mathbb{N}}]$ for some $\left(  y_{n}\right)  \prec\left(
x_{n}\right)  .$ A normalized sequence $(x_{n})$ is said to be an $\ell^{1}%
$-${\mathcal{S}}_{\beta}$-\emph{spreading model with constant} $K$ if
$\Vert\sum_{n\in F}a_{n}x_{n}\Vert\geq K^{-1}\sum_{n\in F}|a_{n}|$ whenever
$F\in{\mathcal{S}}_{\beta}$.

\bigskip

\noindent\textbf{Partly modified mixed Tsirelson spaces}

Let $\left(  \theta_{n}\right)  $ be a null sequence in the interval $(0,1)$
and $\sigma_{n}\in\left\{  U,M\right\}  $ for every $n.$ We say that a family
$\left(  E_{i}\right)  _{i=1}^{k}$ of subsets of \thinspace$\mathbb{N}$ is
$\left(  {\mathcal{S}}_{n},\sigma_{n}\right)  $-adapted if $\left(  \min
E_{i}\right)  _{i=1}^{k}\in{\mathcal{S}}_{n}$ and
\[
\left\{
\begin{array}
[c]{ccc}%
E_{i}\cap E_{j}=\emptyset, & 1\leq i\neq j\leq k & \text{if }\sigma_{n}=M,\\
E_{i}<E_{i+1}, & 1\leq i<k & \text{if }\sigma_{n}=U.
\end{array}
\right.
\]
An $\left(  {\mathcal{S}}_{n},\sigma_{n}\right)  $-adapted family $\left(
E_{i}\right)  _{i=1}^{k}$ is said to be ${\mathcal{S}}_{n}$-\emph{admissible}
(respectively ${\mathcal{S}}_{n}$-\emph{allowable}) if $\sigma_{n}=U$
(respectively $\sigma_{n}=M$). Define the \emph{partly modified mixed
Tsirelson space} $X=T\left[  \left(  {\mathcal{S}}_{n},\sigma_{n},\theta
_{n}\right)  _{n=1}^{\infty}\right]  $ to be the completion of $c_{00}$ under
the implicitly defined norm%
\begin{equation}
\left\Vert x\right\Vert =\max\left\{  \left\Vert x\right\Vert _{c_{0}}%
,\sup_{n}\theta_{n}\sup\sum_{i}\left\Vert E_{i}x\right\Vert \right\}  ,
\label{equation 0.1}%
\end{equation}
where the last supremum is taken over all $\left(  {\mathcal{S}}_{n}%
,\sigma_{n}\right)  $-adapted families $\left(  E_{i}\right)  .$ If
$\sigma_{n}=U$ for all $n$ (respectively $\sigma_{n}=M$ for all $n$), then $X$
is a mixed Tsirelson space (respectively modified mixed Tsirelson space). We
will assume that $\sigma_{p_{0}}=M$ for some $p_{0}.$

\bigskip

\noindent\textbf{Norming Trees}

Equation (\ref{equation 0.1}) can be viewed as an iterative prescription for
computing the norm. The procedure may be summarized in terms of norming trees,
from which the existence and uniqueness of a norm satisfying equation
(\ref{equation 0.1}) also follows. An ($(\mathcal{S}_{n},\sigma_{n})_{n}%
$-)\emph{adapted tree} $\mathcal{T}$ is a finite collection of elements
$(E_{i}^{m})$, $0\leq m\leq r,$ $1\leq i\leq k(m)$, in $[{\mathbb{N}%
}]^{<\infty}$ with the following properties.

\begin{enumerate}
\item $k(0) = 1$,

\item Every $E_{i}^{m+1}$ is a subset of some $E_{j}^{m}$,

\item For each $j$ and $m$, the collection $\{E_{i}^{m+1}:E_{i}^{m+1}\subseteq
E_{j}^{m}\}$ is $\left(  {\mathcal{S}}_{k},\sigma_{k}\right)  $-adapted for
some $k$.
\end{enumerate}

The set $E_{1}^{0}$ is called the \emph{root} of the adapted tree. The
elements $E_{i}^{m}$ are called \emph{nodes} of the tree. If $E_{i}%
^{n}\subseteq E_{j}^{m}$ and $n>m$, we say that $E_{i}^{n}$ is a
\emph{descendant} of $E_{j}^{m}$ and $E_{j}^{m}$ is an \emph{ancestor} of
$E_{i}^{n}$. If, in the above notation, $n=m+1$, then $E_{i}^{n}$ is said to
be an \emph{immediate successor} of $E_{j}^{m}$, and $E_{j}^{m}$ the
\emph{immediate predecessor }or\emph{ parent} of $E_{i}^{n}$. Nodes with no
descendants are called \emph{terminal nodes} or \emph{leaves} of the tree. The
collection of all leaves of $\mathcal{T}$ is denoted by $\mathcal{L}\left(
\mathcal{T}\right)  $. Assign \emph{tags} to the individual nodes inductively
as follows. Let $t(E_{1}^{0})=1$. If $t(E_{i}^{m})$ has been defined and the
collection $(E_{j}^{m+1})$ of all immediate successors of $E_{i}^{m}$ forms an
$\left(  {\mathcal{S}}_{k},\sigma_{k}\right)  $-adapted collection, then
define $t(E_{j}^{m+1})=\theta_{k}t(E_{i}^{m})$ for all immediate successors
$E_{j}^{m+1}$ of $E_{i}^{m}.$ If $x\in c_{00}$ and ${\mathcal{T}}$ is an
adapted tree, let ${\mathcal{T}}x=\sum t(E)\Vert Ex\Vert_{c_{0}}$ where the
sum is taken over all leaves in ${\mathcal{T}}$. It follows from the implicit
description (equation (\ref{equation 0.1})) of the norm in $X$ that $\Vert
x\Vert=\max{\mathcal{T}}x$, with the maximum taken over the set of all adapted
trees. Let us also point out that if ${\mathcal{E}}$ is a collection of
pairwise disjoint nodes of an adapted tree ${\mathcal{T}}$ so that
$E\subseteq\cup{\mathcal{E}}$ for every leaf $E$ of ${\mathcal{T}}$ and $x\in
c_{00}$, then ${\mathcal{T}}x\leq\sum_{F\in{\mathcal{E}}}t(F)\Vert Fx\Vert$.
Given a node $E\in\mathcal{T}$ with tag $t\left(  E\right)  =\prod_{i=1}%
^{m}\theta_{n_{i}},$ define $\operatorname{ord}_{\mathcal{T}}\left(  E\right)
=\sum_{i=1}^{m}n_{i}$. When there is no confusion, we write
$\operatorname{ord}(E)$ instead of $\operatorname{ord}_{\mathcal{T}}(E).$

Let $\mathcal{T}$ be an adapted tree. A node $E\in\mathcal{T}$ is said to be a
\emph{sibling} of $F\in\mathcal{T}$ if they have the same parent. If $(z_{i})$
is a block sequence, we say that $E$ \emph{begins at} $z_{k}$ if
$E\cap\operatorname*{supp}z_{k}\neq\emptyset$ and $E\cap\operatorname*{supp}%
z_{j}=\emptyset$ for all $j<k.$ To say that $E$ \emph{begins before} $z_{k}$
means that $E$ begins at $z_{j}$ for some $j<k$ and we denote this condition
by $E\vartriangleleft z_{k}.$

\bigskip

\noindent$\ell^{1}$\textbf{-Trees and Bourgain's }$\ell^{1}$-\textbf{Index}

A \emph{tree} in a Banach space $B$ is a subset ${\mathcal{T}}$ of $\cup
_{n=1}^{\infty}B^{n}$ so that $(x_{1},\dots,x_{n})\in{\mathcal{T}}$ whenever
$(x_{1},\dots,x_{n},x_{n+1})\in{\mathcal{T}}$. Elements of the tree are called
\emph{nodes}. It is \emph{well-founded} if there is no infinite sequence
$(x_{n})$ so that $(x_{1},\dots,x_{m})\in{\mathcal{T}}$ for all $m$. If $B$
has a basis, then a tree ${\mathcal{T}}$ is said to be a \emph{block tree}
(with respect to the basis) if every node is a block basis of the given basis.
For any well-founded tree ${\mathcal{T}}$, its \emph{derived tree} is the tree
${\mathcal{D}}^{(1)}({\mathcal{T}})$ consisting of all nodes $(x_{1}%
,\dots,x_{n})$ so that $(x_{1},\dots,x_{n},x)\in{\mathcal{T}}$ for some $x$.
Inductively, set ${\mathcal{D}}^{(\alpha+1)}({\mathcal{T}})={\mathcal{D}%
}^{(1)}({\mathcal{D}}^{(\alpha)}({\mathcal{T}}))$ for all ordinals $\alpha$
and ${\mathcal{D}}^{(\alpha)}({\mathcal{T}})=\cap_{\beta<\alpha}{\mathcal{D}%
}^{(\beta)}({\mathcal{T}})$ for all limit ordinals $\alpha$. The \emph{order}
of a tree ${\mathcal{T}}$ is the smallest ordinal $\operatorname*{o}%
({\mathcal{T}})=\alpha$ such that ${\mathcal{D}}^{(\alpha)}({\mathcal{T}%
})=\emptyset$.\newline

\noindent\textbf{Definition}.\ Given a finite constant $K$ $\geq1$, an
$\ell^{1}$-$K$-tree in a Banach space $B$ is a tree in $B$ so that every node
$(x_{1},\dots,x_{n})$ is a normalized sequence such that $\Vert\sum a_{k}%
x_{k}\Vert\geq K^{-1}\sum|a_{k}|$ for all $(a_{k})$. If $B$ has a basis, an
$\ell^{1}$-$K$-block tree is a block tree that is also an $\ell^{1}$-$K$-tree.
Suppose that $B$ does not contain $\ell^{1}$, let $I(B,K)=\sup
\operatorname*{o}({\mathcal{T}})$, where the sup is taken over the set of all
$\ell^{1}$-$K$-trees in $X$. The \emph{Bourgain} $\ell^{1}$-\emph{index} of
$B$ is defined to be $I(B)=\sup_{K<\infty}I(B,K)$. The \emph{block} $\ell^{1}%
$-\emph{index} $I_{b}(B)$ is defined analogously using block trees if $B$ has
a basis. In \cite[Lemmas 5.7 and 5.11]{JO} , it was shown that $I_{b}(B)\ne
I_{b}(B,K)$ and $I(B)\ne I(B,K)$ for every $K$. In particular, $I_{b}(B),
I(B)$ are limit ordinals. It was also shown that \cite[Corollary 5.13]{JO}
$I(B)=I_{b}(B)$ when both are defined and either one has value $\geq
\omega^{\omega}$.

\section{Sufficient conditions for subsequential minimality}

The purpose of the present section is to give sufficient conditions for a
partly modified mixed Tsirelson space to be subsequentially minimal. Prior
experience with mixed Tsirelson spaces \cite{LT1} informs us that there may be
some connection with the Bourgain $\ell^{1}$-index. This indeed turns out to
be the case but the proof requires a different approach.

The main result of the section is the following theorem. The smallest integer
greater than or equal to $a \in{\mathbb{R}}$ is denoted by $\lceil a \rceil$.
For the rest of the section, $X$ will denote a partly modified mixed Tsirelson space.

\begin{thm}
\label{Equivalence} Let $X$ be a partly modified mixed Tsirelson space. If
$Y\prec X$ and $I(Y)>\omega^{\omega},$ then there exists $\left(
x_{n}\right)  \prec Y$ such that $\left(  x_{n}\right)  \sim\left(  e_{p_{n}%
}\right)  $, where $p_{n}=\min\operatorname*{supp}x_{n}.$ Consequently, $X$ is
subsequentially minimal if $I(Y) > \omega^{\omega}$ for all $Y \prec X$.
\end{thm}

Before proceeding with the proof of the theorem, let us draw the following corollary.

\begin{cor}
\label{cor 2} Suppose that there exists $\varepsilon>0$ such that
\[
\sup\{\frac{n}{m}: \theta_{n} \geq\varepsilon^{m}\} = \infty.
\]
Then $X$ is subsequentially minimal. This holds in particular if $\sup
\theta_{n}^{1/n} = 1$.
\end{cor}

\begin{proof}
Clearly, for any $n \in{\mathbb{N}}$ and any $Y \prec X$, every normalized
block sequence in $Y$ is an $\ell^{1}$-${\mathcal{S}}_{n}$-spreading model
with constant $\theta_{n}^{-1}$. By \cite{JO}, if $Y$ contains an $\ell^{1}%
$-${\mathcal{S}}_{2n}$-spreading model with constant $K$, then it contains an
$\ell^{1}$-${\mathcal{S}}_{n}$-spreading model with constant $\sqrt{K}$. With
the assumption of the corollary, for any $k \in{\mathbb{N}}$, there are $m, n$
so that $n/m \geq2k$ and $\theta_{n} \geq\varepsilon^{m}$. Choose $i$ and $j$
so that $2^{i} \leq m < 2^{i+1}$ and $2^{j} \leq n < 2^{j+1}$. Then any $Y
\prec X$ contains an $\ell^{1}$-${\mathcal{S}}_{2^{j}}$-spreading model with
constant $\theta_{n}^{-1}$, and hence, by the remark above, an $\ell^{1}%
$-${\mathcal{S}}_{2^{j-i}}$-spreading model with constant $\theta
_{n}^{-1/2^{i}}$. Since $\theta_{n}^{-1/2^{i}} \leq\varepsilon^{-2}$ and
$2^{j-i} \geq k$, $Y$ has an $\ell^{1}$-${\mathcal{S}}_{k}$-spreading model
with constant $\varepsilon^{-2}$ for all $k$. Hence there is an $\ell^{1}%
$-$\varepsilon^{-2}$-tree on $Y$ of order $\omega^{\omega}$. Thus
$I_{b}(Y,\varepsilon^{-2})\geq\omega^{\omega}$ and so $I(Y) = I_{b}%
(Y)>I_{b}(Y,\varepsilon^{-2})\geq\omega^{\omega}$. The desired result now
follows from Theorem \ref{Equivalence}.

Finally, assume that $\sup\theta_{n}^{1/n} = 1$. Given $0 < \varepsilon< 1$
and $k \in{\mathbb{N}}$, there exists $n > k$ such that $\theta_{n}^{1/n} >
\varepsilon^{1/k}$. Set $m = \lceil n/k\rceil\geq2$. Then $\theta_{n}
\geq\varepsilon^{m}$ and $n/m \geq k(1- 1/m) \geq k/2$.
\end{proof}

The proof of Theorem \ref{Equivalence} occurs in two stages. First we show
that from any block subspace of $X$ with a high $\ell^{1}$-index a
``slow-growing" block sequence may be extracted (see property ($*$) defined
below). In the second part, we show that this subsequence is equivalent to a
subsequence of the unit vector basis $(e_{k})$.\newline

\noindent\textbf{Definition}.\ Let $Y=\left[  \left(  y_{k}\right)  \right]  $
be a block subspace of $X$, we say that $Y$ has property $\left(  \ast\right)
$ if there exists a constant $C<\infty$ such that for all $n\in{\mathbb{N}}$,
there exists a normalized vector $x\in Y_{n}=\left[  \left(  y_{k}\right)
_{k=n}^{\infty}\right]  $ such that $\sum\Vert E_{i}x\Vert\leq C$ whenever
$(E_{i})$ is ${\mathcal{S}}_{n}$-allowable.\newline

First we recall a needed lemma.

\begin{lem}
[{\cite[Proposition 14]{LT}}]\label{lemma 2} Let ${\mathcal{T}}$ be a
well-founded block tree in a Banach space $B$ with a basis. Define
\[
{\mathcal{H}}=\{(\max\operatorname{supp}x_{j})_{j=1}^{r}:(x_{j})_{j=1}^{r}%
\in{{\mathcal{T}}}\}
\]
and
\[
{\mathcal{G}}=\{G:\text{$G$ is a spreading of a subset of some $H\in
{\mathcal{H}}$}\}.
\]
Then ${\mathcal{G}}$ is hereditary and spreading. If ${\mathcal{G}}$ is
compact, then $\iota({\mathcal{G}})\geq$ $\operatorname{o}({{\mathcal{T}}})$.
\end{lem}

\begin{lem}
\label{lemma 3} If $I(Y)>\omega^{\omega}$ then $Y$ has property $\left(
\ast\right)  $.
\end{lem}

\begin{proof}
There exists $K<\infty$ such that $I_{b}(Y,K)\geq\omega^{\omega}$. Let
${\mathcal{T}}$ be an $\ell^{1}$-$K$-block tree in $Y$ such that
$\operatorname*{o}({\mathcal{T}})\geq\omega^{\omega}$. Given $n\geq p_{0}$,
consider the tree $\widehat{{\mathcal{T}}}$ consisting of all nodes of the
form $(x_{j})_{j=n}^{r}$ for some $(x_{j})_{j=1}^{r}\in{\mathcal{T}}$, $r\geq
n$. Then $\widehat{{\mathcal{T}}}$ is an $\ell^{1}$-$K$-block tree in $Y_{n}$
such that $\operatorname*{o}(\widehat{{\mathcal{T}}})\geq\omega^{\omega}$.
Define
\[
{\mathcal{H}}=\{(\max\operatorname{supp}x_{j})_{j=n}^{r}:(x_{j})_{j=n}^{r}%
\in\widehat{{\mathcal{T}}}\}
\]
and
\[
{\mathcal{G}}=\{G:\text{$G$ is a spreading of a subset of some $H\in
{\mathcal{H}}$}\}.
\]
By Lemma \ref{lemma 2}, ${\mathcal{G}}$ is hereditary and spreading, and
either ${\mathcal{G}}$ is noncompact or it is compact with $\iota
({\mathcal{G}})\geq$ $\operatorname*{o}(\widehat{{\mathcal{T}}})\geq
{\omega^{\omega}>\omega}^{n+1}$. By \cite[Theorem 1.1]{G}, there exists
$M\in\lbrack{\mathbb{N}}]$ such that
\[
{\mathcal{S}}_{n+1}\cap\lbrack M]^{<\infty}\subseteq{\mathcal{G}}.
\]
Now \cite[Proposition 3.6]{OTW} gives a finite set $G\in{\mathcal{S}}%
_{n+1}\cap\lbrack M]^{<\infty}$ and a sequence of positive numbers
$(a_{p})_{p\in G}$ such that $\sum a_{p}=1$ and $\sum_{p\in F}a_{p}<\left(
\theta_{p_{0}}\right)  ^{P},$ where $P=\left\lceil \frac{n}{p_{0}}\right\rceil
$, whenever $F\subseteq G$ and $F\in{\mathcal{S}}_{n}$. By definition, there
exist a node $(x_{j})_{j=n}^{r}\in\widehat{{\mathcal{T}}}$ and a subset $J$ of
the integer interval $[n,r]$ such that $G$ is a spreading of $(\max
\operatorname{supp}x_{j})_{j\in J}$. Denote the unique order preserving
bijection from $J$ onto $G$ by $u$ and consider the vector $y=\sum_{j\in
J}a_{u(j)}x_{j}$. Since $(x_{j})_{j=n}^{r}$ is a normalized $\ell^{1}$%
-$K$-block sequence in $Y_{n}$ and $\sum a_{u(j)}=1$, $y\in Y_{n}$ and $\Vert
y\Vert\geq1/K$.

Let $(E_{i})$ be ${\mathcal{S}}_{n}$-allowable. Let $J_{1}=\left\{  j\in
J:\text{some }E_{i}\text{ begins at }x_{j}\right\}  $ and $J_{2}=J\setminus
J_{1}.$ Note that ${\mathcal{S}}_{n}\subseteq{\mathcal{S}}_{p_{0}P}=\left[
{\mathcal{S}}_{p_{0}}\right]  ^{P}$. Thus for each $j$, $\left(  \theta
_{p_{0}}\right)  ^{P}\sum_{i}\left\Vert E_{i}x_{j}\right\Vert \leq\left\Vert
x_{j}\right\Vert =1.$ Also, since $\left\{  u\left(  j\right)  :j\in
J_{1}\right\}  \in{\mathcal{S}}_{n},$ $\sum_{j\in J_{1}}a_{u\left(  j\right)
}<\left(  \theta_{p_{0}}\right)  ^{P}.$ Hence
\begin{align}
\sum_{i}\left\Vert E_{i}\sum_{j\in J_{1}}a_{u\left(  j\right)  }%
x_{j}\right\Vert  &  \leq\sum_{j\in J_{1}}a_{u\left(  j\right)  }\sum
_{i}\left\Vert E_{i}x_{j}\right\Vert \label{J_1}\\
&  \leq\sum_{j\in J_{1}}a_{u\left(  j\right)  }\frac{1}{\left(  \theta_{p_{0}%
}\right)  ^{P}}<1.\nonumber
\end{align}
On the other hand, the collection $\left\{  E_{i}\cap\operatorname*{supp}%
x_{j}:E_{i}\vartriangleleft x_{j}\right\}  $ of pairwise disjoint sets is
$\mathcal{S}_{1}$-allowable and thus $\mathcal{S}_{p_{0}}$-allowable.
Therefore,%
\begin{align}
\sum_{i}\left\Vert E_{i}\sum_{j\in J_{2}}a_{u\left(  j\right)  }%
x_{j}\right\Vert  &  \leq\sum_{j\in J_{2}}a_{u\left(  j\right)  }\sum
_{i}\left\Vert E_{i}x_{j}\right\Vert \label{J_2}\\
&  =\sum_{j\in J_{2}}a_{u\left(  j\right)  }\sum_{E_{i}\vartriangleleft x_{j}%
}\left\Vert E_{i}x_{j}\right\Vert \nonumber\\
&  =\sum_{j\in J_{2}}a_{u\left(  j\right)  }\sum_{E_{i}\vartriangleleft x_{j}%
}\left\Vert \left(  E_{i}\cap\operatorname*{supp}x_{j}\right)  x_{j}%
\right\Vert \nonumber\\
&  \leq\sum_{j\in J_{2}}a_{u\left(  j\right)  }\frac{1}{\theta_{p_{0}}%
}.\nonumber
\end{align}
Combining inequalities (\ref{J_1}) and (\ref{J_2}) gives%
\begin{align*}
\sum\Vert E_{i}y\Vert &  =\sum\Vert E_{i}\sum_{j\in J}a_{u(j)}x_{j}\Vert\\
&  \leq\sum\Vert E_{i}\sum_{j\in J_{1}}a_{u(j)}x_{j}\Vert+\sum\Vert E_{i}%
\sum_{j\in J_{2}}a_{u(j)}x_{j}\Vert\\
&  \leq1+\frac{1}{\theta_{p_{0}}}.
\end{align*}
It is clear that the normalized element $x=y/\Vert y\Vert$ satisfies the
statement of the lemma with the constant $C=\left(  1+\frac{1}{\theta_{p_{0}}%
}\right)  K$.
\end{proof}

We record the quantitative statement of Lemma \ref{lemma 3} for future reference.

\begin{lem}
\label{quant} Let ${\mathcal{T}}$ be an $\ell^{1}$-$K$-block tree on a block
subspace $Y$ of $X$ of order $\operatorname{o}({\mathcal{T}}) \geq
\omega^{\omega}$. Then for all $n \in{\mathbb{N}}$, there is a normalized
vector $x$ in the span of a node of ${\mathcal{T}}$ such that $\sum\|E_{i}x\|
\leq K(1+\theta_{p_{0}}^{-1})$ whenever $(E_{i})$ is ${\mathcal{S}}_{n}$-allowable.
\end{lem}

For each $n \in{\mathbb{N}}$, define
\[
\xi_{n} = \sup\{\theta_{m_{1}}\cdots\theta_{m_{j}}: m_{1} + \cdots+ m_{j} >
n\}.
\]
Then $(\xi_{n})$ is a null sequence. Assume that $Y$ has property $\left(
\ast\right)  ,$ choose $\left(  x_{k}\right)  \prec Y$ and a strictly
increasing sequence $\left(  n_{k}\right)  ,$ $n_{0}=1,$ so that for each $k,$

\begin{enumerate}
\item[$\left(  \alpha\right)  $] $\sum\Vert E_{s}x_{k}\Vert\leq C$ whenever
$(E_{s})$ is ${\mathcal{S}}_{n_{k-1}}$-allowable$,$

\item[$\left(  \beta\right)  $] $\xi_{n_{k}}\left\Vert x_{k}\right\Vert
_{\ell^{1}}\leq\frac{1}{2^{k}}$,

\item[($\gamma$)] $2q_{k} \leq p_{k+1}$ for all $k$, where $p_{k} =
\min\operatorname{supp} x_{k}$ and $q_{k} = \max\operatorname{supp} x_{k}$.
\end{enumerate}

Let $\left(  b_{k}\right)  _{k=1}^{N}\in c_{00}^{+}$ and set $x=\sum_{k=1}%
^{N}b_{k}x_{k}.$

\begin{lem}
\label{Tree} Let $\mathcal{T}$ be an adapted tree. If $\mathcal{E}$ is a
collection of pairwise disjoint nodes of $\mathcal{T}$ such that
$\operatorname{ord}(E)\leq m$ for all $E\in\mathcal{E}$, then $\mathcal{E}$ is
$\mathcal{S}_{m}$-allowable.
\end{lem}

\begin{proof}
Note that if $\mathcal{T}$ is an adapted tree, then it is an allowable tree
with nodes of the same orders. The conclusion follows from \cite[Lemma 3]{LT2}.
\end{proof}

\begin{lem}
\label{Trimming}Given any adapted tree $\mathcal{T}$, there exists an adapted
tree $\mathcal{T}^{\prime}$ such that

\begin{enumerate}
\item[(a)] if $E\in\mathcal{T}^{\prime}$ and $E\cap\operatorname*{supp}%
x_{k}\neq\emptyset,$ then $\operatorname{ord}\left(  E\right)  \leq n_{k},$

\item[(b)] $\mathcal{T}x\leq\mathcal{T}^{\prime}x+\sum\frac{b_{k}}{2^{k}}.$
\end{enumerate}
\end{lem}

\begin{proof}
Given an adapted tree $\mathcal{T}$ and $F\subseteq\mathbb{N}$, define
\[
\mathcal{T}_{F}=\left\{  E\cap F:E\in\mathcal{T}\text{, }E\cap F\neq
\emptyset\right\}  .
\]
Clearly $\mathcal{T}_{F}$ is an adapted tree. For all $k=2,...,N,$, define a
set $F_{k}$ by
\[
F_{k}^{c}=\cup\left\{  E\cap\operatorname*{supp}x_{k}:E\in\mathcal{T},
\operatorname{ord}\left(  E\right)  >n_{k}\right\}  .
\]
Then
\begin{align*}
\mathcal{T}x_{k}  &  =\sum_{E\in\mathcal{L}\left(  \mathcal{T}\right)
}t\left(  E\right)  \left\Vert Ex_{k}\right\Vert \\
&  =\sum_{\substack{E\in\mathcal{L}\left(  \mathcal{T}\right)
\\\operatorname{ord}\left(  E\right)  \leq n_{k}}}t\left(  E\right)
\left\Vert Ex_{k}\right\Vert +\sum_{\substack{E\in\mathcal{L}\left(
\mathcal{T}\right)  \\\operatorname{ord}\left(  E\right)  >n_{k}}}t\left(
E\right)  \left\Vert Ex_{k}\right\Vert \\
&  \leq\mathcal{T}_{F_{k}}x_{k}+\xi_{n_{k}}\left\Vert x_{k}\right\Vert
_{\ell^{1}}\leq\mathcal{T}_{F_{k}}x_{k}+\frac{1}{2^{k}}.
\end{align*}
Let $\mathcal{T}^{\prime}=\mathcal{T}_{\left(  F_{2}\cap F_{3}\cap\dots\cap
F_{N}\right)  }.$ Note that $\mathcal{T}^{\prime}$ satisfies (a) and
$\mathcal{T}^{\prime}x_{k}=\mathcal{T}_{F_{k}}x_{k}$ if $2\leq k\leq N.$ Hence%
\begin{align*}
\mathcal{T}x  &  =\sum b_{k}\mathcal{T}x_{k}\\
&  \leq b_{1}\mathcal{T}x_{1}+\sum_{k=2}^{N}b_{k}\left(  \mathcal{T}_{F_{k}%
}x_{k}+\frac{1}{2^{k}}\right) \\
&  =\sum b_{k}\mathcal{T}^{\prime}x_{k}+\sum\frac{b_{k}}{2^{k}}\\
&  =\mathcal{T}^{\prime}x+\sum\frac{b_{k}}{2^{k}}.
\end{align*}

\end{proof}

Define $\mathcal{E}_{k}$ to be the set
\[
\{E\in\mathcal{T}^{\prime}: E \text{ begins at $x_{k}$ and has a sibling that
begins before } x_{k}\}.
\]

\begin{lem}
\label{Estimate}$\sum\limits_{E\in\mathcal{E}_{k}}\left\Vert Ex_{k}\right\Vert
\leq C$\ for all $k=2,...,N$.
\end{lem}

\begin{proof}
Note that if $E\in\mathcal{E}_{k},$ $E$ has a sibling $E^{\prime}$ that begins
before $x_{k}.$ Hence $\operatorname{ord}\left(  E\right)  =
\operatorname{ord}\left(  E^{\prime}\right)  \leq n_{k-1}$ by property (a) of
Lemma \ref{Trimming}. By Lemma \ref{Tree}, $\mathcal{E}_{k}$ is $\mathcal{S}%
_{n_{k-1}}$-allowable. The conclusion follows from condition $\left(
\alpha\right)  .$
\end{proof}

\begin{proof}
[Proof of Theorem \ref{Equivalence}]As $(e_{k})$ is a $1$-unconditional basis
of $X$, it is enough to consider nonnegative coefficients. As above, consider
$\left(  b_{k}\right)  _{k=1}^{N}\in c_{00}^{+}$ and set $x=\sum_{k=1}%
^{N}b_{k}x_{k}$, $y = \sum_{k=1}^{N}b_{k}e_{p_{k}}$. It is easy to see that
$\|y\| \leq\|x\|$. We will show that $\|x\| \leq(2+C)\|y\|$, where $C$ is the
constant in condition ($\alpha$). Given an adapted tree ${\mathcal{T}}$, we
obtain an adapted tree ${\mathcal{T}}^{\prime}$ as in Lemma \ref{Trimming}. We
may further assume that every node $E \in{\mathcal{T}}^{\prime}\backslash
\mathcal{L}({\mathcal{T}}^{\prime})$ is the union of its immediate successors,
that $E \subseteq\cup_{k}\operatorname{supp} x_{k}$ for every $E
\in{\mathcal{T}}^{\prime}$ and that, relabeling if necessary, the root of
${\mathcal{T}}^{\prime}$ begins at $x_{1}$. With these assumptions, every node
$E \in\mathcal{L}({\mathcal{T}}^{\prime})$ that intersects
$\operatorname{supp} x_{k}$, $k \geq2$, is a descendant of some node in
$\mathcal{E}_{k}$.
For each $k\geq2,$ choose $E_{k}\in\mathcal{E}_{k}$ such that $t\left(
E_{k}\right)  =\max\left\{  t\left(  E\right)  :E\in\mathcal{E}_{k}\right\}
.$ By Lemma \ref{Estimate}, for $k \geq2$,
\[
\mathcal{T}^{\prime}x_{k} =\sum_{E\in\mathcal{E}_{k}}t\left(  E\right)
\left\Vert Ex_{k}\right\Vert \leq t\left(  E_{k}\right)  \sum_{E\in
\mathcal{E}_{k}}\left\Vert Ex_{k}\right\Vert \leq t\left(  E_{k}\right)  C.
\]
Therefore,
\begin{align*}
\mathcal{T}x  &  \leq\mathcal{T}\left(  b_{1}x_{1}\right)  +\sum_{k=2}%
^{N}b_{k}\mathcal{T}^{\prime}x_{k}+\sum_{k=2}^{N}\frac{b_{k}}{2^{k}}\\
&  \leq b_{1}+C\sum_{k=2}^{N}b_{k}t\left(  E_{k}\right)  +\sum_{k=2}^{N}%
\frac{b_{k}}{2^{k}}\\
&  \leq C\sum_{k=2}^{N}t\left(  E_{k}\right)  b_{k}+2\left\Vert \left(
b_{k}\right)  \right\Vert _{c_{0}}\\
&  \leq C\sum_{k=2}^{N}t\left(  E_{k}\right)  b_{k}+2\left\Vert y\right\Vert .
\end{align*}
To complete the proof, it suffices to appeal to Proposition \ref{tree} below
to see that $\sum_{k=2}^{N}t\left(  E_{k}\right)  b_{k} \leq\|y\|$.
\end{proof}

\noindent\textbf{Remark.} This proof above shows that if $(x_{k})$ is a
(possibly finite) normalized block sequence in $X$ satisfying conditions
($\alpha$), ($\beta$) and ($\gamma$) for some $(n_{k})$, then $(x_{k})$ is
$(2+C)$-equivalent to $(e_{p_{k}})$.

\begin{prop}
\label{tree} There is an $\left(  \mathcal{S}_{n},\sigma_{n}\right)
_{n=1}^{\infty}$-adapted tree $\mathcal{T}^{\prime\prime}$ so that
\[
\mathcal{T}^{\prime\prime}y\geq\sum_{k=1}^{N}b_{k}t(E_{k}).
\]
In particular, $\sum_{k=1}^{N}b_{k}t(E_{k}) \leq\|y\|$.
\end{prop}

\noindent The tree $\mathcal{T}^{\prime\prime}$ is constructed by substituting
each node $E$ in $\mathcal{T}^{\prime}$ with one or two nodes, which we now
proceed to describe. For each $E\in\mathcal{T}^{\prime},$ define
$G_{E}=\left\{  p_{j}:E_{j}\subsetneqq E\right\}  $. If $E \in{\mathcal{T}%
}^{\prime}$ and $E \neq E_{k}$ for any $k$, substitute $G_{E}$ for $E$. If $E
= E_{k}$ for some $k$, substitute two nodes, namely $\{p_{k}\}$ and $G_{E}$,
in place of $E$. The resulting collection of nodes after the substitutions we
denote by ${\mathcal{T}}^{\prime\prime}$. Note that since the root of
${\mathcal{T}}^{\prime}$ begins at $x_{1}$, it cannot be equal to $E_{k}$ for
any $k$. Thus the root of ${\mathcal{T}}^{\prime}$ is substituted with a
single node. To show that $\mathcal{T}^{\prime\prime}$ \ is an $\left(
\mathcal{S}_{n},\sigma_{n}\right)  _{n=1}^{\infty}$-adapted tree, it is enough
to show that if $E\in\mathcal{T}^{\prime}$ has immediate successors $\left(
F_{i}\right)  _{i=1}^{s}$ which form an $\left(  \mathcal{S}_{n},\sigma
_{n}\right)  $-adapted family$,$ then $\left(  G_{F_{i}}\right)  ^{s}_{i=1}
\cup P$ is an $\left(  \mathcal{S}_{n},\sigma_{n}\right)  $-adapted family of
subsets of $G_{E},$ where $P=\{\left\{  p_{k}\right\}  :F_{i}=E_{k}\text{ for
some }i\}$. We divide the proof of this assertion into a series of claims and
lemmas.\newline

\noindent\underline{Claim 1.} $\left(  G_{F_{i}}\right)  ^{s}_{i=1} \cup P$ is
a family of pairwise disjoint subsets of $G_{E}.$

By definition, $\{p_{k}\}\subseteq G_{E}$ for any $\{p_{k}\} \in P$. Let us
show that $G_{F_{i}}\subseteq G_{E}.$ Indeed, if $p_{j}\in G_{F_{i}},$ then
$E_{j}\subsetneqq F_{i}\subseteq E$. Thus $p_{j}\in G_{E}.$

Now if $i\neq i^{\prime},$ then $F_{i}\cap F_{i^{\prime}}=\emptyset$. By
definition, $G_{F_{i}}$ is disjoint from $G_{F_{i^{\prime}}}$. If $F_{i}%
=E_{k}$ for some $i$ and $k$, then for any $i^{\prime}$ (including $i$
itself), $E_{k}\subsetneqq F_{i^{\prime}}$ cannot hold. Therefore, $\{p_{k}\}$
and $G_{F_{i^{\prime}}}$ are disjoint for all $i^{\prime}$. Since obviously
any two sets in $P$ are disjoint, the claim is established.\newline

\noindent\underline{Claim 2.} If $(F_{i})^{s}_{i=1}$ consists of successive
sets, then so does $\left(  G_{F_{i}}\right)  ^{s}_{i=1} \cup P$.

First we show that if $F_{i} < F_{i^{\prime}}$, then $G_{F_{i}} <
G_{F_{i^{\prime}}}$. Let $p_{j}\in G_{F_{i}}$ and $p_{j^{\prime}}\in
G_{F_{i^{\prime}}}$. Then $E_{j}\subsetneqq F_{i}$ and $E_{j^{\prime}%
}\subsetneqq F_{i^{\prime}}.$ Since $E_{j}$ begins at $x_{j}$, $E_{j^{\prime}%
}$ begins at $x_{j^{\prime}}$ and $F_{i}<F_{i^{\prime}},$ it follows that $j <
j^{\prime}$ and hence $p_{j}<p_{j^{\prime}}.$ This shows that $G_{F_{i}%
}<G_{F_{i^{\prime}}}.$

Next, if $F_{i} < F_{i^{\prime}}=E_{k}$ for some $i, i^{\prime}$ and $k$, then
we claim that $G_{F_{i}}<\left\{  p_{k}\right\}  <G_{F_{i^{\prime}}}.$ To see
the first inequality, pick a point $p_{j}\in G_{F_{i}}$. Then $E_{j}%
\subsetneqq F_{i}.$ In particular, $E_{j}<F_{i^{\prime}}=E_{k}.$ Since $E_{j}$
begins at $x_{j}$ and $E_{k}$ begins at $x_{k},$ we deduce that $j<k$ and thus
$p_{j}<p_{k}$. Hence $G_{F_{i}} <\left\{  p_{k}\right\}  .$ Similarly, if
$p_{j}\in G_{F_{i^{\prime}}},$ then $E_{j}\subsetneqq F_{i^{\prime}}=E_{k}.$
Since $E_{j}$ begins at $x_{j}$ and $E_{k}$ begins at $x_{k},$ we deduce that
$k<j.$ This shows that $\left\{  p_{k}\right\}  <G_{F_{i^{\prime}}}.$\newline

Let $\hat{P} = \{p_{k}: \{p_{k}\} \in P\}$.\newline

\noindent\underline{Claim 3.} $\left(  \min G_{F_{i}}\right)  _{i=1}^{s}%
\cup\hat{P} \in\mathcal{S}_{n}$

The proof of this claim requires several short lemmas.

\begin{lem}
\label{Lemma 10} For any $E\in\mathcal{T}^{\prime}$, $\min G_{E}\geq2\min E.$
\end{lem}

\begin{proof}
Suppose that $p_{j}\in G_{E}.$ Then $E_{j}\subsetneqq E.$ Since $E_{j}$ has a
sibling that begins before $x_{j},$ $E$ begins before $x_{j}.$ This implies
that%
\[
2\min E \leq2q_{j-1} \leq p_{j}\quad\text{ by }(\gamma).
\]

\end{proof}

\begin{lem}
\label{Lemma 11} $\hat{P}$ is a spreading of a subset of $(\min F_{i})
_{i=1}^{s}$ and $p_{k}\geq2\min F_{1}$ for all $p_{k}\in\hat{P}.$
\end{lem}

\begin{proof}
We may assume that $\min F_{1}<\cdots<\min F_{s}.$ For each $k,$ let
$H_{k}=\left\{  \min F_{i}:\min F_{i}\in\operatorname*{supp}x_{k}\right\}  .$
List the $k$'s such that $H_{k}\neq\emptyset$ in increasing order as
$k_{1}<\cdots<k_{r}.$ Since every $F_{i}$ begins at or after $x_{k_{1}},$
$E_{k_{1}}\neq F_{i}$ for any $i.$ Therefore, $\hat{P}\subseteq\left(
p_{k_{\ell}}\right)  _{\ell=2}^{r}$. For each $2\leq\ell\leq r,$ choose
$i_{\ell-1}$ such that $\min F_{i_{\ell-1}}\in H_{k_{\ell-1}}.$ Then
$(p_{k_{\ell}}) _{\ell=2}^{r}$ is a spreading of $(\min F_{i_{\ell-1}})
_{\ell=2}^{r}.$ Also note that $p_{k}\geq p_{k_{2}}\geq2 q_{k_{1}}\geq2\min
F_{1}$ for all $p_{k}\in\hat{P}.$
\end{proof}

It follows from Lemmas \ref{Lemma 10} and \ref{Lemma 11} that $(\min G_{F_{i}%
}) _{i=1}^{s}\cup\hat{P}$ can be written as $\cup_{j\in B}A_{j},$ where $B =
\{2\min F_{1}\} \cup(\min F_{i})^{s}_{i=2}$, $\min A_{j}\geq j$, and $\vert
A_{j}\vert\leq2$ for all $j\in B.$

\begin{lem}
\label{union} Suppose that $n\in\mathbb{N}$, $L\in\mathcal{S}_{n}$ and $B$ is
a spreading of $L$ such that $\min B\geq2\min L.$ If $\left\vert
A_{j}\right\vert \leq2$ and $\min A_{j}\geq j$ for all $j\in B,$ then
$\cup_{j\in B}A_{j}\in\mathcal{S}_{n}.$
\end{lem}

\begin{proof}
It is easy to see that we may assume $A_{j}<A_{j^{\prime}}$ if $j<j^{\prime}$.
Write $L=\cup_{k=1}^{p}L_{k},$ where $L_{1}<\cdots<L_{p}$ are in
$\mathcal{S}_{n-1}$ and $p\leq\min L_{1}.$ Then $B=\cup_{k=1}^{p}B_{k},$ where
each $B_{k}$ is a spreading of $L_{k}$ and $B_{1}<\cdots<B_{p}.$ Denoting by
$\mathcal{A}_{2}$ the collection of subsets of ${\mathbb{N}}$ having at most
two elements, we appeal to \cite[Remark on p.312]{LT} to deduce that
\[
\cup_{j\in B_{k}}A_{j}\in\mathcal{S}_{n-1}\left[  \mathcal{A}_{2}\right]
\subseteq\left(  \mathcal{S}_{n-1}\right)  ^{2}.
\]
Hence $\cup_{j\in B}A_{j}=\cup_{k=1}^{2p}C_{i},$ where $C_{1}<\cdots<C_{2p}$
are in $\mathcal{S}_{n-1}.$ Since $2p\leq2\min L_{1}\leq\min B\leq\min C_{1},$
the conclusion of the lemma follows.
\end{proof}

\begin{proof}
[Completion of proof of Proposition \ref{tree}]It follows from the claims and
lemmas above that the nodes of ${\mathcal{T}}^{\prime\prime}$ form an
$(S_{n},\sigma_{n})^{\infty}_{n=1}$-adapted tree, where the tag of any node in
${\mathcal{T}}^{\prime\prime}$ is the same as the tag of the node in
${\mathcal{T}}^{\prime}$ for which it is a substitute. Moreover, it follows
from Claim 1 that all nodes in $P$ are terminal. Therefore,
\[
{\mathcal{T}}^{\prime\prime}y \geq\sum_{\{p_{k}\}\in P}t(\{p_{k}\})b_{k} =
\sum^{N}_{k=2}t(E_{k})b_{k}.
\]

\end{proof}

Recall that a Banach space $Z$ is said to be \emph{minimal} if every infinite
dimensional subspace of $Z$ has a further subspace isomorphic to $Z$. This
definition is due to Rosenthal. In \cite{Go}, Gowers introduced the more
general notion of \emph{quasi-minimal} spaces. Two Banach spaces are said to
be \emph{totally incomparable} if they do not have isomorphic infinite
dimensional subspaces. A Banach space is said to be quasi-minimal if it does
not contain a pair of totally incomparable infinite dimensional closed
subspaces. Using Theorem \ref{Equivalence}, Corollary \ref{cor 2} and
Proposition \ref{qm} below, we obtain

\begin{cor}
Let $X = T[({\mathcal{S}}_{n},\sigma_{n},\theta_{n})_{n=1}^{\infty}] $ be a
partly modified mixed Tsirelson space so that $I(Y)>\omega^{\omega}$ for every
block subspace $Y$ of $X$. Then $X$ is quasi-minimal. This holds if there
exists $\varepsilon>0$ such that $\sup\{n/m: \theta_{n} \geq\varepsilon^{m}\}
= \infty$, and, in particular, if $\sup\theta_{n}^{1/n} = 1$.
\end{cor}

\begin{prop}
\label{qm} Let $(p_{k})$ and $(q_{k})$ be subsequences of ${\mathbb{N}}$ so
that $p_{k} \leq q_{k} < 2q_{k} \leq p_{k+1}$ for all $k$. Then the sequences
$(e_{p_{k}})$ and $(e_{q_{k}})$ are $2$-equivalent in any partly modified
mixed Tsirelson space $X = T[({\mathcal{S}}_{n},\sigma_{n},\theta_{n}%
)_{n=1}^{\infty}]$.
\end{prop}

\begin{proof}
Define a sequence of norms on $X$ follows. Let $\Vert x\Vert_{0}=\Vert
x\Vert_{c_{0}}$ and
\[
\Vert x\Vert_{i+1}=\max\{\Vert x\Vert_{0},\sup_{n}\sup\theta_{n}\sum_{m}\Vert
E_{m}x\Vert_{i}\},
\]
where the final supremum is taken over all $({\mathcal{S}}_{n},\sigma_{n}%
)$-adapted families $(E_{m})$. It is clear that $\Vert x\Vert=\lim\Vert
x\Vert_{i}$ for all $x\in X$. For any finite subset $E$ of $(q_{k})$, let the
\emph{shift} of $E$ be the set $s(E)=\{p_{k}:q_{k}\in E\}$. We claim that for
any $i$, any $(a_{k})\in c_{00}$ and any $E\subseteq(q_{k})$, there exist
$p_{j}\in s(E)$ and $F\subseteq s(E)$ such that $p_{j}<F$ and%
\begin{equation}
\Vert E\sum a_{k}e_{q_{k}}\Vert_{i}\leq|a_{j}|+\Vert F\sum a_{k}e_{p_{k}}%
\Vert_{i}.\label{claim}%
\end{equation}
Once the claim is proved, it follows easily that $\Vert\sum a_{k}e_{q_{k}%
}\Vert\leq2\Vert\sum a_{k}e_{p_{k}}\Vert$. Since each ${\mathcal{S}}_{n}$ is
spreading, we clearly have $\Vert\sum a_{k}e_{p_{k}}\Vert\leq\Vert\sum
a_{k}e_{q_{k}}\Vert$, and the proof of the proposition would be complete. We
now prove the claim (\ref{claim}) by induction on $i$. The case $i=0$ is
trivial. Suppose that the claim holds for some $i$. We may assume that
\[
\Vert E\sum a_{k}e_{q_{k}}\Vert_{i+1}=\theta_{n}\sum_{m=1}^{d}\Vert E_{m}\sum
a_{k}e_{q_{k}}\Vert_{i},
\]
where $(E_{m})_{m=1}^{d}$ is an $({\mathcal{S}}_{n},\sigma_{n})$-adapted
family of subsets of $E$, arranged so that $(\min E_{m})_{m=1}^{d}$ is an
increasing sequence. By induction, for each $m$, there are $p_{j_{m}}\in
s(E_{m})$ and $F_{m}\subseteq s(E_{m})$ such that $p_{j_{m}}<F_{m}$ and
\[
\Vert E_{m}\sum a_{k}e_{q_{k}}\Vert_{i}\leq|a_{j_{m}}|+\Vert F_{m}\sum
a_{k}e_{p_{k}}\Vert_{i}.
\]
Observe that for every $m$, $2\min E_{m}\leq2q_{j_{m}}<p_{j_{m}+1}\leq\min
F_{m}$. Also, for $m\geq2$, $2\min E_{m-1}\leq\min s(E_{m})\leq p_{j_{m}}$.
Let $m_{0}$ be such that $p_{j_{m_{0}}}$ is the minimum of the sequence
$(p_{j_{m}})_{m=1}^{d}$. Then $(p_{j_{m}})_{m\neq m_{0}}\cup(\min F_{m}%
)_{m=1}^{d}$ may be written as $\cup_{j\in B}A_{j}$, where $B$ is a spreading
of $(\min E_{m})_{m=1}^{d}$ such that $\min B\geq2\min E_{1}$, $|A_{j}|\leq2$
and $A_{j}\geq j$ for all $j\in B$. By Lemma \ref{union}, $(p_{j_{m}})_{m\neq
m_{0}}\cup(\min F_{m})_{m=1}^{d}\in{\mathcal{S}}_{n}$. Clearly, $\{\{p_{j_{m}%
}\}:m\neq m_{0}\}\cup\{F_{m}:1\leq m\leq d\}$ is a pairwise disjoint family
that is successive if $(E_{m})_{m=1}^{d}$ is. Thus, this family is
$({\mathcal{S}}_{n},\sigma_{n})$-adapted. We may then conclude that
\begin{align*}
\Vert E\sum a_{k}e_{q_{k}}\Vert_{i+1} &  =\theta_{n}\sum_{m=1}^{d}\Vert
E_{m}\sum a_{k}e_{q_{k}}\Vert_{i}\\
&  \leq\theta_{n}|a_{j_{m_{0}}}|+\theta_{n}(\sum_{m\neq m_{0}}|a_{j_{m}}%
|+\sum_{m=1}^{d}\Vert F_{m}\sum a_{k}e_{p_{k}}\Vert_{i})\\
&  \leq|a_{j_{m_{0}}}|+\Vert F\sum a_{k}e_{p_{k}}\Vert_{i+1},
\end{align*}
where $F=\{p_{j_{m}}:m\neq m_{0}\}\cup\cup_{m=1}^{d}F_{m}\subseteq s(E)$ and
$F>p_{j_{m_{0}}}\in s(E)$.
\end{proof}

If $X=T[({\mathcal{S}}_{n},\sigma_{n},\theta_{n})_{n=1}^{\infty}]$ is a partly
modified mixed Tsirelson space where $\sigma_{p_{0}}=M$, then it is clear that
every disjointly supported sequence $(x_{k})_{k=1}^{n}$ in $[(e_{k}%
)_{k=n}^{\infty}]$ is $\theta_{p_{0}}^{-1}$-equivalent to the unit vector
basis of $\ell^{1}(n)$. Such spaces are called \emph{strongly asymptotic
$\ell^{1}$ spaces}. In \cite{DFKO}, it was proved that every minimal, strongly
asymptotic $\ell^{1}$ Banach space with a basis is isomorphic to a subspace
$\ell^{1}$. Since partly modified spaces are reflexive (this may be proved
using the arguments of \cite{ADKM}; alternatively, it follows from the
computation of the $\ell^{1}$-index below (Theorem \ref{UppBdIndex})), we get
that no partly modified mixed Tsirelson space contains a minimal subspace.
Hence the class of the partly modified mixed Tsirelson spaces $X$ such that
$I(Y)>\omega^{\omega}$ for every subspace $Y$ of $X$ provides examples of
quasi-minimal Banach spaces without minimal subspaces.

\section{The Bourgain $\ell^{1}$-index}

In this section, we develop the techniques in \S 2 further to investigate the
Bourgain $\ell^{1}$-index of partly modified mixed Tsirelson spaces. In the
first part of the section, we show that $I\left(  X\right)  $ does not exceed
$\omega^{\omega\cdot2}$. In the second part, we pinpoint the value of
$I\left(  X\right)  $ in certain cases in terms of the sequence of
coefficients $\left(  \theta_{n}\right)  $.

In the following proposition, we will require the concepts of block subtrees,
minimal trees ${\mathcal{T}}_{\alpha}$ and replacement trees ${\mathcal{T}%
}(\alpha,\beta)$ defined, constructed and developed in \cite{JO}. We refer the
reader to that paper for details. The execution of the following proof is
comparable to that of \cite[Lemma 4.2]{JO}. When two trees $\mathcal{T}$ and
$\mathcal{T}^{\prime}$ are isomorphic, we write $\mathcal{T}\simeq
\mathcal{T}^{\prime}$. Given two finite sequences $\vec{x}=\left(
x_{1},\cdots,x_{m}\right)  $ and $\vec{y}=\left(  y_{1},\cdots,y_{n}\right)
,$ let $\vec{x}\sqcup\vec{y}=\left(  x_{1},\cdots,x_{m},y_{1},\cdots
,y_{n}\right)  .$ We say that a normalized vector $x$ satisfies property
$(\ast)$ for the couple $(n,C)\in\mathbb{N\times R}^{+}$ if $\sum\left\Vert
E_{i}x\right\Vert \leq C$ whenever $(E_{i})$ is $\mathcal{S}_{n}$-allowable.

\begin{prop}
\label{Replacement}If $\mathcal{T}$ is an $\ell^{1}$-$K$-block tree of order
$\operatorname*{o}({\mathcal{T}})\geq\omega^{\omega}\cdot\alpha,$ then for any
$n_{0}\in\mathbb{N}$ and any positive sequence $(\varepsilon_{i})$, there
exists a block subtree ${\mathcal{T}}^{\prime}$ of ${\mathcal{T}}$, isomorphic
to ${\mathcal{T}}_{\alpha}$, such that every node $\left(  x_{1},\cdots
,x_{d}\right)  \in{\mathcal{T}}^{\prime}$ satisfies

\begin{enumerate}
\item There exist $n_{1}<\cdots<n_{d-1}$, with $n_{1}>n_{0}$, such that each
$x_{i}$ satisfies property $\left(  \ast\right)  $ for the couple $\left(
n_{i-1},C\right)  ,$ where $C=\left(  1+\theta_{p_{0}}^{-1}\right)  K$,

\item $\xi_{n_{i}}\left\Vert x_{i}\right\Vert _{\ell^{1}}\leq\varepsilon_{i}$
for $1\leq i<d$, and

\item $2\max\operatorname{supp} x_{i} \leq\min\operatorname{supp} x_{i+1}$ if
$1 \leq i < d$.
\end{enumerate}
\end{prop}

\begin{proof}
The proof is by induction on $\alpha.$ The case $\alpha=1$ follows from Lemma
\ref{quant}. Suppose that $\mathcal{T}$ is an $\ell^{1}$-$K$-block tree of
order $\operatorname*{o}({\mathcal{T}})\geq\omega^{\omega}\cdot\left(
\alpha+1\right)  .$ According to \cite[Lemma 3.7]{JO} and replacing
${\mathcal{T}}$ by a subtree if necessary, we may assume that $\mathcal{T}$ is
isomorphic to the \textquotedblleft replacement tree" $\mathcal{T}$ $\left(
\alpha+1,\omega^{\omega}\right)  .$ From the definition of $\mathcal{T}\left(
\alpha+1,\omega^{\omega}\right)  ,$ we see that $\left(  \mathcal{T}\left(
\alpha+1,\omega^{\omega}\right)  \right)  ^{\left(  \omega^{\omega}\cdot
\alpha\right)  }$ is the minimal tree $\mathcal{T}_{\omega^{\omega}}$.
Applying the case $\alpha=1$ to $\mathcal{T}^{\left(  \omega^{\omega}%
\cdot\alpha\right)  }\simeq\mathcal{T}_{\omega^{\omega}}$, we obtain a
normalized block $y$ of a node $\vec{x}=\left(  x_{1},\cdots,x_{m}\right)  $
in $\mathcal{T}^{\left(  \omega^{\omega}\cdot\alpha\right)  }$ such that $y$
satisfies ($\ast$) for the couple $\left(  n_{0},C\right)  .$ Choose
$n_{1}>n_{0}$ such that $\xi_{n_{1}}\left\Vert y\right\Vert _{\ell^{1}}%
\leq\varepsilon_{1}.$ Without loss of generality, we may assume that $\vec{x}$
is a terminal node in $\mathcal{T}^{\left(  \omega^{\omega}\cdot\alpha\right)
}.$ By the construction of $\mathcal{T}\left(  \alpha+1,\omega^{\omega
}\right)  ,$ the subtree $\mathcal{T}_{\vec{x}}$ of $\mathcal{T}$ consisting
of all nodes $\vec{z}>\vec{x}$ is isomorphic to $\mathcal{T}\left(
\alpha,\omega^{\omega}\right)  $ and hence has order $\omega^{\omega}%
\cdot\alpha.$ Consider the \textquotedblleft restricted subtree"
$\mathcal{R}\left(  \mathcal{T}_{\vec{x}}\right)  $ \cite[Definition 4.1]{JO}
consisting of all $(w_{j}, \dots, w_{k})$, where $\vec{x}\sqcup(w_{1}%
,\dots,w_{k})\in\mathcal{T}_{\vec{x}}$ and $j$ is the smallest integer such
that $\min\operatorname{supp} w_{j} \geq2 \max\operatorname{supp} x_{m}$. Then
$\mathcal{R}\left(  \mathcal{T}_{\vec{x}}\right)  $ is an $\ell^{1}$-$K$-block
tree of order $\omega^{\omega}\cdot\alpha.$ Apply the inductive hypothesis to
$\mathcal{R}\left(  \mathcal{T}_{\vec{x}}\right)  $ with the parameters
$n_{1}$ and $(\varepsilon_{i+1})$ to obtain a block subtree $\mathcal{T}%
^{\prime\prime}$ of $\mathcal{R}\left(  \mathcal{T}_{\vec{x}}\right)  .$
Define $\mathcal{T}^{\prime}=\left\{  (y)\sqcup\vec{w}:\vec{w}\in
\mathcal{T}^{\prime\prime}\right\}  .$ It is easy to check that $\mathcal{T}%
^{\prime}$ satisfies the desired conclusion (for the ordinal $\alpha+1$).

Suppose that $\mathcal{T}$ is an $\ell^{1}$-$K$-block tree of order
$\operatorname*{o}({\mathcal{T}})\geq\omega^{\omega}\cdot\alpha,$ where
$\alpha$ is a limit ordinal. Let $\left(  \alpha_{n}\right)  $ be a sequence
of ordinals strictly increasing to $\alpha.$ Then $\mathcal{T}$ contains
pairwise disjoint subtrees $\mathcal{T}_{n}$ with $\operatorname*{o}%
({\mathcal{T}}_{n})\geq\omega^{\omega}\cdot\alpha_{n}$ for all $n.$ For each
$n$, apply the inductive hypothesis to obtain a block subtree $\mathcal{T}%
_{n}^{\prime}$ of $\mathcal{T}_{n}.$ The block subtree $\mathcal{T}^{\prime
}=\cup\mathcal{T}_{n}^{\prime}$ of $\mathcal{T}$ satisfies the conclusion of
the proposition.
\end{proof}

If $(\varepsilon_{i})$ is chosen to be $(1/2^{i})$, then from the remark
following the proof of Theorem \ref{Equivalence}, we see that every node
$\left(  x_{1},\cdots,x_{d}\right)  \in{\mathcal{T}}^{\prime}$ is
$(2+C)$-equivalent to $\left(  e_{p_{i}}\right)  $, where $p_{i}%
=\min\operatorname*{supp}x_{i}.$ For $y \in c_{00}$, let $\|y\|_{\mathcal{S}%
_{p}} = \sup_{E\in\mathcal{S}_{p}}\|Ey\|_{\ell^{1}}$.

\begin{thm}
\label{UppBdIndex}The Bourgain $\ell^{1}$-index of $X=T\left[  \left(
{\mathcal{S}}_{n},\sigma_{n},\theta_{n}\right)  _{n=1}^{\infty}\right]  $ is
$I\left(  X\right)  \leq\omega^{\omega\cdot2}.$
\end{thm}

\begin{proof}
If $I\left(  X\right)  >\omega^{\omega\cdot2},$ then by \cite[Corollary
5.13]{JO}, there exists an $\ell^{1}$-$K$-block tree $\mathcal{T}$ with
$\operatorname*{o}({\mathcal{T}})\geq\omega^{\omega\cdot2}$ for some $K>0.$
Let $n$ be chosen so that $\xi_{n}<\frac{1}{2K\left(  2+C\right)  }.$ By
Proposition \ref{Replacement}, we obtain an $\ell^{1}$-$K$-block tree
${\mathcal{T}}^{\prime}$ of ${\mathcal{T}}$ with $\operatorname*{o}%
({\mathcal{T}}^{\prime})=\omega^{n+1}$ such that every node $\left(
x_{1},\cdots,x_{d}\right)  $ in ${\mathcal{T}}^{\prime}$ is $\left(
2+C\right)  $-equivalent to $\left(  e_{p_{i}}\right)  $. Define
\[
{\mathcal{H}}=\{(p_{j})_{j=n}^{r}:(x_{j})_{j=n}^{r}\in{\mathcal{T}}^{\prime
}\}
\]
and
\[
{\mathcal{G}}=\{G:\text{$G$ is a spreading of a subset of some $H\in
{\mathcal{H}}$}\}.
\]
By Lemma \ref{lemma 2}, ${\mathcal{G}}$ is hereditary and spreading, and
either ${\mathcal{G}}$ is noncompact or it is compact with $\iota
({\mathcal{G}})\geq$ $\operatorname*{o}({\mathcal{T}}^{\prime})\geq
{\omega^{n+1}>\omega}^{n}$. By \cite[Theorem 1.1]{G}, there exists
$M\in\lbrack{\mathbb{N}}]$ such that ${\mathcal{S}}_{n}\cap\lbrack
M]^{<\infty}\subseteq{\mathcal{G}}$. As in the proof of Lemma \ref{lemma 3},
we obtain a node $(x_{j})^{r}_{j=n} \in{\mathcal{T}}^{\prime}$, $J
\subseteq[n,r]$, an order preserving map $u$ from $J$ onto a spreading of
$(p_{j})_{j\in J}$ and a sequence of positive numbers $(a_{u(j)})_{j\in J}$
such that $\sum_{j\in J}a_{u\left(  j\right)  }=1$ and $\sum_{j\in
A}a_{u\left(  j\right)  }<{\xi_{n}}$ whenever $\left\{  u\left(  j\right)
:j\in A\right\}  \in\mathcal{S}_{n-1}.$ Let $y = \sum_{j\in J}a_{u(j)}x_{j}$.
Since $(x_{j})_{j}$ is a normalized $\ell^{1}$-$K$-block sequence, $\Vert
y\Vert\geq1/K.$ On the other hand,
\begin{align*}
\Vert y\Vert &  =\Vert\sum_{j\in J}a_{u(j)}x_{j}\Vert\leq\left(  2+C\right)
\Vert\sum_{j\in J}a_{u(j)}e_{p_{i}}\Vert\\
&  \leq\left(  2+C\right)  (\|\sum_{j\in J}a_{u(j)}e_{p_{j}}\|_{{\mathcal{S}%
}_{n-1}}+ \xi_{n}\Vert(a_{u(j)})\Vert_{\ell^{1}})\leq2\left(  2+C\right)
\xi_{n},
\end{align*}
contradicting the choice of $n.$
\end{proof}

In the second half of the section, we obtain an estimate on the norms of
vectors spanned by normalized block sequences in $X$ (Proposition
\ref{Prop B4}), from which the value of the Bourgain $\ell^{1}$-index $I(X)$
may be deduced. For the remainder of the section, assume that $\left(
x_{k}\right)  $ is a normalized block sequence in $X = T[{\mathcal{S}}_{n},
\sigma_{n}, \theta_{n})^{\infty}_{n=1}]$, $\left(  a_{k}\right)  \in c_{00}$
and $q_{k}=\max\operatorname*{supp}x_{k}.$ Set $x=\sum a_{k}x_{k}.$ Recall the
assumption that $\sigma_{p_{0}}=M$ for some $p_{0}.$ Given a node $E$ in an
adapted tree ${\mathcal{T}}$, we say that it is a \emph{long node} (with
respect to $x$) if $E\cap\operatorname*{supp}x_{k}\neq\emptyset$ for more than
one $k$. Otherwise, we term the node \emph{short}.

\begin{lem}
\label{Lemma B1}For any $N$, there exists an adapted tree $\mathcal{T}$ such
that all long nodes $E\in\mathcal{T}$ satisfy $t\left(  E\right)  >\theta_{N}$
and
\[
\left\Vert x\right\Vert \leq\mathcal{T}x+\frac{\theta_{N}}{\theta_{p_{0}}%
}\left\Vert \left(  a_{k}\right)  \right\Vert _{\ell^{1}}.
\]

\end{lem}

\begin{proof}
Choose an adapted tree $\mathcal{T}^{\prime}$ such that $\left\Vert
x\right\Vert =\mathcal{T}^{\prime}x.$ Let $\mathcal{E}$ be the collection of
minimal elements in the set of long nodes $E$ with $t\left(  E\right)
\leq\theta_{N}$. For each $E\in\mathcal{E},$ let $k_{E}$ be the smallest $k$
such that $\operatorname*{supp}x_{k}\cap E\neq\emptyset\ $\ and let
$F_{E}=\operatorname*{supp}x_{k_{E}}\cap E.$ For each $k,$ the nonempty sets
in the collection $\left\{  \left(  E\smallsetminus F_{E}\right)
\cap\operatorname*{supp}x_{k}\right\}  $ is $\mathcal{S}_{1}$-allowable and
hence $\mathcal{S}_{p_{0}}$-allowable. Thus,%
\[
\sum_{E\in\mathcal{E}}t\left(  E\right)  \left\Vert \left(  E\smallsetminus
F_{E}\right)  x_{k}\right\Vert \leq\theta_{N}\sum_{E\in\mathcal{E}}\left\Vert
\left(  E\smallsetminus F_{E}\right)  x_{k}\right\Vert \leq\frac{\theta_{N}%
}{\theta_{p_{0}}}.
\]
Then%
\[
\sum_{E\in\mathcal{E}}t\left(  E\right)  \left\Vert \left(  E\smallsetminus
F_{E}\right)  x\right\Vert \leq\frac{\theta_{N}}{\theta_{p_{0}}}\left\Vert
\left(  a_{k}\right)  \right\Vert _{\ell^{1}}.
\]
Let $\mathcal{T}$ be the tree obtained from $\mathcal{T}^{\prime}$ by changing
all nodes $G\in\mathcal{T}^{\prime}$, $G\subseteq E$ for some $E\in
\mathcal{E}$ to $G\cap F_{E},$ Then $\mathcal{T}$ is an adapted tree such that
all long nodes $H$ in $\mathcal{T}$ satisfies $t\left(  H\right)  >\theta
_{N}.$ Moreover,%
\begin{align*}
\left\Vert x\right\Vert  &  =\mathcal{T}^{\prime}x\leq\mathcal{T}x+\sum
_{E\in\mathcal{E}}t\left(  E\right)  \left\Vert \left(  E\smallsetminus
F_{E}\right)  x\right\Vert \\
&  \leq\mathcal{T}x+\frac{\theta_{N}}{\theta_{p_{0}}}\left\Vert \left(
a_{k}\right)  \right\Vert _{\ell^{1}}.
\end{align*}

\end{proof}

Fix $N$ and let ${\mathcal{T}}$ be the tree given by Lemma \ref{Lemma B1}. For
any $\varepsilon>0,$ let $k\left(  \varepsilon\right)  =\max\left\{
n_{1}+\cdots+n_{j}:\theta_{n_{1}}\cdots\theta_{n_{j}}>\varepsilon\right\}  .$
Let $\mathcal{E}$ denote the set of all minimal short nodes in $\mathcal{T}$.

\begin{lem}
\label{Lemma B2}If $\mathcal{E}_{1}=\left\{  E\in\mathcal{E}:E\text{ has a
long sibling}\right\}  ,$ then
\[
\sum_{E\in\mathcal{E}_{1}}t\left(  E\right)  \left\Vert Ex\right\Vert
\leq\left\Vert \sum a_{k}e_{q_{k}}\right\Vert _{\mathcal{S}_{k\left(
\theta_{N}\right)  }}.
\]

\end{lem}

\begin{proof}
If $E\in\mathcal{E}_{1},$ then $t\left(  E\right)  >\theta_{N}$ and hence
$\operatorname{ord}\left(  E\right)  \leq k\left(  \theta_{N}\right)  .$ Hence
by Lemma \ref{Tree},
$\mathcal{E}_{1}$ is $\mathcal{S}_{k\left(  \theta_{N}\right)  }$-allowable.
Since each $E\in\mathcal{E}_{1}$ is a short node, it follows that the set
$Q_{0}=\left\{  q_{k}:\operatorname*{supp}x_{k}\cap E\neq\emptyset\text{ for
some }E\in\mathcal{E}_{1}\right\}  \in\mathcal{S}_{k\left(  \theta_{N}\right)
}.$ Thus
\[
\sum_{E\in\mathcal{E}_{1}}t\left(  E\right)  \left\Vert Ex\right\Vert \leq
\sum_{q_{k}\in Q_{0}}\left\vert a_{k}\right\vert \leq\left\Vert \sum
a_{k}e_{q_{k}}\right\Vert _{\mathcal{S}_{k\left(  \theta_{N}\right)  }}.
\]

\end{proof}

For $m,n\in\mathbb{N}$, define $\eta_{m,n}=\inf\frac{\theta_{m+n}}%
{\theta_{n_{1}}\cdots\theta_{n_{s}}},$ where the infimum is taken over all
$n_{1},\dots,n_{s}$ such that $n_{1}+\cdots+n_{s}\geq n,$ with the additional
requirement that $\sigma_{n_{1}}=\sigma_{n_{2}}=\cdots=\sigma_{n_{s}}=M$ if
$\sigma_{m+n}=M.$

\begin{lem}
\label{Lemma B3}Suppose that $\inf\limits_{m}\limsup\limits_{n}\eta_{m,n}=0.$
For any $\varepsilon>0,$ there exist $m$ and $n_{0}$ such that
\[
\sum_{E\in\mathcal{E}\backslash\mathcal{E}_{1}}t\left(  E\right)  \left\Vert
Ex\right\Vert \leq\varepsilon\left\Vert \left(  a_{k}\right)  \right\Vert
_{\ell^{1}}+2\left\Vert \sum a_{k}e_{q_{k}}\right\Vert _{\mathcal{S}_{k\left(
\theta_{N}\right)  }+n_{0}+m}.
\]

\end{lem}

\begin{proof}
Choose $m$ and $n_{0}$ such that $\eta_{m,n}<\varepsilon$ if $n\geq n_{0}.$
Let $\mathcal{D}=\left(  D_{i}\right)  $ be the set of all parents of nodes in
$\mathcal{E\smallsetminus E}_{1}.$ In particular, each $D_{i}$ is a long node
and hence $t\left(  D_{i}\right)  >\theta_{N}.$ It follows that
$\operatorname{ord}\left(  D_{i}\right)  \leq k\left(  \theta_{N}\right)  $.
Also, the nodes in $\mathcal{D}$ are pairwise disjoint since no $E\in
\mathcal{E\smallsetminus E}_{1}$ can have a long sibling. For each $i$, there
exists some $n_{i}$ such that $\mathcal{F}_{i}=\left\{  E\in
\mathcal{E\smallsetminus E}_{1}:E\subseteq D_{i}\right\}  $ is $\left(
\mathcal{S}_{n_{i}},\sigma_{n_{i}}\right)  $-adapted. Let $I=\left\{
i:n_{i}\leq n_{0}+m\right\}  .$ Then $\operatorname{ord}\left(  E\right)
=\operatorname{ord}\left(  D_{i}\right)  +n_{i}\leq k\left(  \theta
_{N}\right)  +n_{0}+m$ for all $E\in{\bigcup\limits_{i\in I}}\mathcal{F}_{i}$.
By Lemma \ref{Tree},
$%
{\displaystyle\bigcup\limits_{i\in I}}
\mathcal{F}_{i}$ is an $\mathcal{S}_{k\left(  \theta_{N}\right)  +n_{0}+m}%
$-allowable collection of short nodes. It follows that%
\[
Q_{0}=\left\{  q_{k}:\operatorname*{supp}x_{k}\cap E\neq\emptyset\text{ for
some }E\in%
{\displaystyle\bigcup\limits_{i\in I}}
\mathcal{F}_{i}\right\}  \in\mathcal{S}_{k\left(  \theta_{N}\right)  +n_{0}%
+m}.
\]
Therefore,%
\begin{equation}
\sum_{E\in\cup_{i\in I}\mathcal{F}_{i}}t\left(  E\right)  \left\Vert
Ex\right\Vert \leq\sum_{q_{k}\in Q_{0}}\left\vert a_{k}\right\vert
\leq\left\Vert \sum a_{k}e_{q_{k}}\right\Vert _{\mathcal{S}_{k\left(
\theta_{N}\right)  +n_{0}+m}}. \label{I}%
\end{equation}
Now consider those $i\notin I.$ Let $\mathcal{F}_{ik}=\left\{  E\in
\mathcal{F}_{i}:E\subseteq\operatorname*{supp}x_{k}\right\}  $. For each $k,$
let
\begin{align*}
I_{k}  &  =\left\{  i\notin I:\left\{  \min E:E\in\mathcal{F}_{ik}\right\}
\in\mathcal{S}_{n_{i}-m}\right\} \\
\text{and}\quad I_{k}^{\prime}  &  =\left\{  i\notin I:\left\{  \min
E:E\in\mathcal{F}_{ik}\right\}  \notin\mathcal{S}_{n_{i}-m}\right\}  .
\end{align*}
Suppose that $i\in I_{k}.$ Choose $n_{1}^{\left(  i\right)  },\dots
,n_{s}^{\left(  i\right)  }\,$such that $n_{1}^{\left(  i\right)  }%
+\cdots+n_{s}^{\left(  i\right)  }\geq n_{i}-m,$
\[
\frac{\theta_{m+n_{i}-m}}{\theta_{n_{1}^{\left(  i\right)  }}\cdots
\theta_{n_{s}^{\left(  i\right)  }}}<\varepsilon
\]
and $\sigma_{n_{1}^{\left(  i\right)  }}=\cdots=\sigma_{n_{s}^{\left(
i\right)  }}=M$ if $\sigma_{n_{i}}=M.$ This is possible since $i\notin I$
implies that $n_{i}-m\geq n_{0}$ and hence $\eta_{m,n_{i}-m}<\varepsilon.$

If $\sigma_{n_{i}}=U,$ then the sets in $\mathcal{F}_{i}$ and hence
$\mathcal{F}_{ik}$ are successive. Since $\left\{  \min E:E\in\mathcal{F}%
_{ik}\right\}  \in\mathcal{S}_{n_{i}-m},$ $\mathcal{F}_{ik}$ is $\mathcal{S}%
_{n_{i}-m}$-admissible and hence $\mathcal{S}_{n_{1}^{\left(  i\right)
}+\cdots+n_{s}^{\left(  i\right)  }}$-admissible. Then%
\begin{equation}
\sum_{E\in\mathcal{F}_{ik}}\theta_{n_{1}^{\left(  i\right)  }}\cdots
\theta_{n_{s}^{\left(  i\right)  }}\left\Vert Ex_{k}\right\Vert =\sum
_{E\in\mathcal{F}_{ik}}\theta_{n_{1}^{\left(  i\right)  }}\cdots\theta
_{n_{s}^{\left(  i\right)  }}\left\Vert ED_{i}x_{k}\right\Vert \leq\left\Vert
D_{i}x_{k}\right\Vert . \label{II}%
\end{equation}
If $\sigma_{n_{i}}=M,$ then $\mathcal{F}_{ik}$ is $\mathcal{S}_{n_{i}-m}%
$-allowable and hence $\mathcal{S}_{n_{1}^{\left(  i\right)  }+\cdots
+n_{s}^{\left(  i\right)  }}$-allowable. Since $\sigma_{n_{1}^{\left(
i\right)  }}=\cdots=\sigma_{n_{s}^{\left(  i\right)  }}=M$, we obtain the same
inequality as (\ref{II}).

From inequality (\ref{II}),
\begin{align*}
\sum_{i\in I_{k}}\sum_{E\in\mathcal{F}_{ik}}t\left(  E\right)  \left\Vert
Ex_{k}\right\Vert  &  =\sum_{i\in I_{k}}t\left(  D_{i}\right)  \theta_{n_{i}%
}\sum_{E\in\mathcal{F}_{ik}}\left\Vert Ex_{k}\right\Vert \\
&  \leq\varepsilon\sum_{i\in I_{k}}t\left(  D_{i}\right)  \theta
_{n_{1}^{\left(  i\right)  }}\cdots\theta_{n_{s}^{\left(  i\right)  }}%
\sum_{E\in\mathcal{F}_{ik}}\left\Vert Ex_{k}\right\Vert \\
&  \leq\varepsilon\sum_{i\in I_{k}}t\left(  D_{i}\right)  \left\Vert
D_{i}x_{k}\right\Vert \leq\varepsilon.
\end{align*}
Therefore,%
\begin{equation}
\sum_{\left\{  \left(  i,k\right)  :i\in I_{k}\right\}  }\sum_{E\in
{\mathcal{F}}_{ik}}t\left(  E\right)  \left\Vert Ex\right\Vert \leq
\varepsilon\left\vert \left\vert \left(  a_{k}\right)  \right\vert \right\vert
_{\ell^{1}}. \label{III}%
\end{equation}
For each $i\notin I$, set $J_{i}=\left\{  k:i\in I_{k}^{\prime}\right\}  .$
Then $\left\{  \min E:E\in\mathcal{F}_{ik}\right\}  \notin\mathcal{S}%
_{n_{i}-m}$ for each $k\in J_{i}$ but $\cup_{k}\left\{  \min E:E\in
\mathcal{F}_{ik}\right\}  =\left\{  \min E:E\in\mathcal{F}_{i}\right\}
\in\mathcal{S}_{n_{i}}.$ By \cite[Lemma 2]{LT}, $\left(  \min\cup
_{E\in\mathcal{F}_{ik}}E\right)  _{k\in J_{i}}\in\mathcal{S}_{m}.$ Now
$\operatorname{ord}\left(  D_{i}\right)  \leq k\left(  \theta_{N}\right)  $
for all $i$ and $\mathcal{D}$ consist of pairwise disjoint sets. Thus by Lemma
\ref{Tree},
$\mathcal{D}$ is $\mathcal{S}_{k\left(  \theta_{N}\right)  }$-allowable.
Therefore, $\left\{  q_{k}:k\in\cup_{i\notin I}J_{i}\right\}  \in
\mathcal{S}_{k\left(  \theta_{N}\right)  +m}.$ It follows that
\begin{align}
\label{IV}\sum_{\{ (i,k) :i\in I_{k}^{\prime}\} }\sum_{E\in\mathcal{F}_{ik}}t(
E) \Vert Ex\Vert &  \leq\sum_{k\in\cup_{i\notin I}J_{i}}\left\vert
a_{k}\right\vert \\
&  \leq\left\Vert \sum a_{k}e_{q_{k}}\right\Vert _{\mathcal{S}_{k\left(
\theta_{N}\right)  +m}}.\nonumber
\end{align}
Combining $\left(  \ref{I}\right)  ,$ $\left(  \ref{III}\right)  $ and
$\left(  \ref{IV}\right)  $ yields%
\[
\sum_{E\in\mathcal{E}\smallsetminus\mathcal{E}_{1}}t\left(  E\right)
\left\Vert Ex\right\Vert =\varepsilon\|(a_{k})\| _{\ell^{1}}+2\left\Vert \sum
a_{k}e_{q_{k}}\right\Vert _{\mathcal{S}_{k\left(  \theta_{N}\right)  +n_{0}
+m}}.
\]

\end{proof}

From Lemmas \ref{Lemma B1}, \ref{Lemma B2}, and \ref{Lemma B3} we have

\begin{prop}
\label{Prop B4}Suppose that $\inf\limits_{m}\limsup\limits_{n}\eta_{m,n}=0.$
Then given any $\varepsilon>0$ and $N,$ there exist $m$ and $n_{0}$ such that%
\[
\left\Vert x\right\Vert \leq\left(  \varepsilon+\frac{\theta_{N}}%
{\theta_{p_{0}}}\right)  \|(a_{k})\| _{\ell^{1}}+3\left\Vert \sum
a_{k}e_{q_{k}}\right\Vert _{\mathcal{S}_{k\left(  \theta_{N}\right)  +n_{0}%
+m}}.
\]

\end{prop}

\begin{thm}
If $\inf\limits_{m}\limsup\limits_{n}\eta_{m,n}=0,$ then $I\left(  X\right)  =
\omega^{\omega}.$
\end{thm}

\begin{proof}
Since $X$ contains $\ell^{1}$-$S_{n}$-spreading models with constant
$\theta_{n}^{-1}$ for all $n$, it is clear that $I(X) \geq\omega^{\omega}$.
Suppose $I(X) > \omega^{\omega}$. There exist $K>1$ and an $\ell^{1}$%
-$K$-block tree ${\mathcal{T}}$ such that $\operatorname*{o}\left(
{\mathcal{T}}\right)  >\omega^{\omega}.$ Let ${\mathcal{H}}\left(
{\mathcal{T}}\right)  =\{(\max\operatorname{supp}x_{j})_{j=1}^{r}%
:(x_{j})_{j=1}^{r}\in{{\mathcal{T}}}\}$ and ${\mathcal{G}}=\{G:G$ is a
spreading of a subset of some $H\in{\mathcal{H}}$$\}.$ Then $\iota
({\mathcal{G}})\geq$ $\operatorname*{o}({{\mathcal{T}}})>\omega^{\omega}.$
Choose $\varepsilon$ and $N$ such that $\varepsilon+\frac{\theta_{N}}%
{\theta_{p_{0}}}<\frac{1}{2K}$ and let $r=k\left(  \theta_{N}\right)
+n_{0}+m$ where $n_{0}, m$ are such that $\eta_{m,n}<\varepsilon$ if $n\geq
n_{0}$. Applying \cite[Theorem 1.1]{G}, there exists $M\in\left[
\mathbb{N}\right]  $ such that $\mathcal{S}_{\omega}\cap\left[  M\right]
^{<\infty}\subseteq{\mathcal{G}}.$ Hence, it follows from \cite[Proposition
3.6]{OTW} that there exist $G=\left(  t_{i}\right)  \in{\mathcal{G}}$ and
$\left(  a_{i}\right)  \in c_{00}^{+}$ such that $\sum a_{i}=1$ and
$\left\Vert \sum a_{i}e_{t_{i}}\right\Vert _{\mathcal{S}_{r}}<\frac{1}{6K}.$

By definition, there exists a normalized $\ell^{1}$-$K$-block sequence
$\left(  x_{i}\right)  _{1}^{k}$ in $X$ such that $\left(  t_{i}\right)  $ is
a spreading of $\left(  q_{i}\right)  =\left(  \max\operatorname*{supp}%
x_{i}\right)  .$ By Proposition \ref{Prop B4},
\begin{align*}
\frac{1}{K}  &  \leq\left\Vert \sum a_{i}x_{i}\right\Vert \leq\frac{1}%
{2K}\left\vert \left\vert \left(  a_{i}\right)  \right\vert \right\vert
_{\ell^{1}}+3\left\Vert \sum a_{i}e_{q_{i}}\right\Vert _{\mathcal{S}_{r}}\\
&  \leq\frac{1}{2K}+3\left\Vert \sum a_{i}e_{t_{i}}\right\Vert _{\mathcal{S}%
_{r}} <\frac{1}{K},
\end{align*}
a contradiction.
\end{proof}

\begin{thm}
\label{Theorem B6}If $\inf\limits_{m}\limsup\limits_{n}\inf\limits_{n_{1}%
+\cdots+n_{s}\geq n}\frac{\theta_{m+n}}{\theta_{n_{1}}\cdots\theta_{n_{s}}%
}>0,$ then $X$ contains $\ell^{1}$-$\mathcal{S}_{m}$-spreading models with
uniform constant. In particular, $I(X) = \omega^{\omega\cdot2}$.
\end{thm}

The proof is based on the following construction.

\begin{lem}
For any $n\in\mathbb{N}$, $\varepsilon>0$ and $L\in\left[  \mathbb{N}\right]
,$ there exists $x\in c_{00}$ such that
\[
\left\Vert x\right\Vert _{\ell^{1}}=\frac{1}{\theta_{n}}%
,\ \operatorname*{supp}x\in\mathcal{S}_{N+1}\cap[L]^{<\infty} \text{ and }
\left\Vert x\right\Vert _{X}\leq1+\frac{1}{\varepsilon},
\]
where $N=\max\left\{  n_{1}+\cdots+n_{s}:\varepsilon\theta_{n_{1}}\cdots
\theta_{n_{s}}>\theta_{n}\right\}  .$ (We take $\max\emptyset=0$.)
\end{lem}

\begin{proof}
According to \cite[Proposition 3.6]{OTW}, there exists $x\in c_{00}$ such that
$\left\Vert x\right\Vert _{\ell^{1}}=\frac{1}{\theta_{n}},$
$\operatorname*{supp}x\in\mathcal{S}_{N+1}\cap\lbrack L]^{<\infty}$ and
$\left\Vert x\right\Vert _{\mathcal{S}_{N}}\leq1.$ If $\mathcal{T}$ is an
adapted tree, then%
\begin{align*}
\mathcal{T}x  &  =\sum_{\substack{E\in\mathcal{L}\left(  \mathcal{T}\right)
\\\varepsilon t\left(  E\right)  \leq\theta_{n}}}t\left(  E\right)  \left\Vert
Ex\right\Vert _{c_{0}}+\sum_{\substack{E\in\mathcal{L}\left(  \mathcal{T}%
\right)  \\\varepsilon t\left(  E\right)  >\theta_{n}}}t\left(  E\right)
\left\Vert Ex\right\Vert _{c_{0}}\\
&  \leq\frac{\theta_{n}}{\varepsilon}\left\Vert x\right\Vert _{\ell^{1}}%
+\sum_{\substack{E\in\mathcal{L}\left(  \mathcal{T}\right)  \\\varepsilon
t\left(  E\right)  >\theta_{n}}}\left\Vert Ex\right\Vert _{c_{0}}.
\end{align*}
But $\varepsilon t\left(  E\right)  >\theta_{N}$ implies that
$\operatorname{ord}\left(  E\right)  \leq N.$ It follows from Lemma \ref{Tree}
that $\left\{  E\in\mathcal{L}\left(  \mathcal{T}\right)  :\text{ }\varepsilon
t\left(  E\right)  >\theta_{n}\right\}  $ is $\mathcal{S}_{N}$-allowable. Then
$\mathcal{T}x\leq\frac{1}{\varepsilon}+\left\Vert x\right\Vert _{\mathcal{S}%
_{N}}\leq\frac{1}{\varepsilon}+1.$
\end{proof}

\begin{proof}
[Proof of Theorem \ref{Theorem B6}]Let $\varepsilon>0$ be such that
\[
\inf\limits_{m}\limsup\limits_{n}\inf\limits_{n_{1}+\cdots+n_{s}\geq n-m}%
\frac{\theta_{n}}{\theta_{n_{1}}\cdots\theta_{n_{s}}}>\varepsilon.
\]
Given any $m,$ choose $n>m$ such that%
\[
\inf\limits_{n_{1}+\cdots+n_{s}\geq n-m}\frac{\theta_{n}}{\theta_{n_{1}}%
\cdots\theta_{n_{s}}}>\varepsilon.
\]
Then $N=\max\left\{  n_{1}+\cdots+n_{s}:\varepsilon\theta_{n_{1}}\cdots
\theta_{n_{s}}>\theta_{n}\right\}  <n-m.$ Choose a block sequence $\left(
x_{k}\right)  $ such that $\left\Vert x_{k}\right\Vert _{\ell^{1}}=\frac
{1}{\theta_{n}},$ $\operatorname*{supp}x_{k}\in\mathcal{S}_{N+1}$ and
$\left\Vert x_{k}\right\Vert _{X}\leq1+\frac{1}{\varepsilon}$ for all $k.$ Let
$p_{k}=\min\operatorname*{supp}x_{k}$ for all $k.$ If $F\in\mathcal{S}_{m},$
then $\left(  p_{k}\right)  _{k\in F}\in\mathcal{S}_{m}$ and hence $%
{\textstyle\bigcup_{k\in F}}
\operatorname*{supp}x_{k}\in\mathcal{S}_{m+N+1}\subseteq\mathcal{S}_{n}.$ Thus
for any $\left(  a_{k}\right)  \in c_{00},$%
\[
\left\Vert \sum_{k\in F}a_{k}x_{k}\right\Vert \geq\theta_{n}\left\Vert
\sum_{k\in F}a_{k}x_{k}\right\Vert _{\ell^{1}}=\sum_{k\in F}\left\vert
a_{k}\right\vert .
\]
This shows that $\left(  x_{k}/\left\Vert x_{k}\right\Vert \right)  $ is an
$\ell^{1}$-$\mathcal{S}_{m}$-spreading model with constant $1+1/\varepsilon.$

Let $K$ be a fixed constant so that for each $m$, there is a normalized block
sequence $(x_{i}^{m})_{i=1}^{\infty}$ that is an $\ell^{1}$-${\mathcal{S}}%
_{m}$-spreading model with constant $K$. If ${\mathcal{F}}$ is a regular
family, consider the tree ${\mathcal{T}}({\mathcal{F}})$ in $X$ consisting of
all sequences of the form $(x_{i}^{m_{1}})_{i\in I_{1}}\cup\cdots\cup
(x_{i}^{m_{r}})_{i\in I_{r}}$ with $I_{k}\in\mathcal{S}_{m_{k}}$, $1\leq k\leq
r$, $i_{k+1}>\max_{i\in I_{k}}\max\operatorname{supp}x_{i}^{m_{k}}$ for all
$i_{k+1}\in I_{k+1}$, $1\leq k<r$, and $(\min I_{1},\dots,\min I_{r}%
)\in{\mathcal{F}}$. If $(x_{i}^{m_{1}})_{i\in I_{1}}\cup\cdots\cup
(x_{i}^{m_{r}})_{i\in I_{r}}\in{\mathcal{T}}({\mathcal{F}}^{(1)})$, take
$i_{0}=\max_{i\in I_{r}}\max\operatorname{supp}x_{i}^{m_{r}}$. There exists
$j_{0}$ so that $(\min I_{1},\dots\min I_{r},j_{0})\in{\mathcal{F}}$. Then
$(x_{i}^{m_{1}})_{i\in I_{1}}\cup\cdots\cup(x_{i}^{m_{r}})_{i\in I_{r}}%
\cup(x_{i}^{m})_{i\in I}\in{\mathcal{T}}({\mathcal{F}})$ provided
$I\in{\mathcal{S}}_{m}$ and $I>\max\{i_{0},j_{0}\}$. It follows easily that
${\mathcal{T}}({\mathcal{F}}^{(1)})\subseteq{\mathcal{T}}({\mathcal{F}%
})^{(\omega^{\omega})}$. Carrying on inductively, one deduces that
$\operatorname{o}({\mathcal{T}}({\mathcal{S}}_{n}))\geq\omega^{\omega}%
\cdot\omega^{n}$ for all $n$. Finally, note that if $(x_{i}^{m_{1}})_{i\in
I_{1}}\cup\cdots\cup(x_{i}^{m_{r}})_{i\in I_{r}}\in{\mathcal{T}}({\mathcal{S}%
}_{n})$, then for all scalars $(a_{i}^{m})$,
\[
\Vert\sum_{k=1}^{r}\sum_{i\in I_{k}}a_{i}^{m_{k}}x_{i}^{m_{k}}\Vert\geq
\theta_{n}\sum_{k=1}^{r}\Vert\sum_{i\in I_{k}}a_{i}^{m_{k}}x_{i}^{m_{k}}%
\Vert\geq\frac{\theta_{n}}{K}\sum_{k=1}^{r}\sum_{i\in I_{k}}|a_{i}^{m_{k}}|.
\]
Hence ${\mathcal{T}}({\mathcal{S}}_{n})$ is an $\ell^{1}$-$K\theta_{n}^{-1}%
$-tree in $X$ of order at least $\omega^{\omega+n}$. Thus $I(X)\geq
\omega^{\omega\cdot2}$. The reverse inequality holds by Theorem
\ref{UppBdIndex}.
\end{proof}

The foregoing results serve to pinpoint the value of the Bourgain $\ell^{1}%
$-index precisely if $X$ is either ``boundedly modified" or ``(completely) modified".

\begin{cor}
Suppose that there exists $N$ such that $\sigma_{n}=U$ for all $n>N,$ or that
$\sigma_{n} = M$ for all $n$. Then

\begin{enumerate}
\item $I\left(  X\right)  =\omega^{\omega}$ if $\inf\limits_{m}\limsup
\limits_{n}\inf\limits_{n_{1}+\cdots+n_{s}\geq n}\frac{\theta_{m+n}}%
{\theta_{n_{1}}\cdots\theta_{n_{s}}}=0,$

\item $I\left(  X\right)  =\omega^{\omega\cdot2}$ if $\inf\limits_{m}%
\limsup\limits_{n}\inf\limits_{n_{1}+\cdots+n_{s}\geq n}\frac{\theta_{m+n}%
}{\theta_{n_{1}}\cdots\theta_{n_{s}}}>0.$ In this case $X$ has $\ell^{1}%
$-$\mathcal{S}_{m}$-spreading models with uniform constant.
\end{enumerate}
\end{cor}

\section{Mixed Tsirelson spaces that are strongly non-subsequentially minimal}

In the final section, it is shown that a large class of (unmodified) mixed
Tsirelson spaces fails to be subsequentially minimal in a strong sense. We
consider a mixed Tsirelson space $X = T[({\mathcal{S}}_{n},\theta_{n}%
)^{\infty}_{n=1}] = T[({\mathcal{S}}_{n},\sigma_{n},\theta_{n})^{\infty}%
_{n=1}]$, where $\sigma_{n} = U$ for all $n$. In this case, we may assume
without loss of generality that $\left(  \theta_{n}\right)  $ is a regular
sequence, i.e., $\left(  \theta_{n}\right)  $ is a non-increasing null
sequence in $\left(  0,1\right)  $ such that $\theta_{m+n}\geq\theta_{n}%
\theta_{m}$ for all $m,n\in\mathbb{N}$. By \cite[Lemma 4.13]{OTW},
$\theta=\lim_{n}\theta_{n}^{1/n}$ exists and is equal to $\sup\theta_{n}%
^{1/n}$. Also, we let $\varphi_{n}=\theta_{n}/\theta^{n}.$\newline

\noindent\textbf{Definition}.\ We say that a Banach space $X$ with a
normalized basis $\left(  e_{k}\right)  $ is \emph{strongly
non-subsequentially minimal }if for every normalized block basis $\left(
x_{k}\right)  $ of $\left(  e_{k}\right)  ,$ there exists $\left(
y_{k}\right)  \prec\left(  x_{k}\right)  $ such that for all $\left(
z_{k}\right)  \prec\left(  y_{k}\right)  ,$ $\left(  z_{k}\right)  $ is not
equivalent to any subsequence of $\left(  e_{k}\right)  .$\newline

The main result of this section is Theorem \ref{Theorem1} where it is shown
that $X$ is strongly non-subsequentially minimal if $\theta< 1$ and $0 <
\inf\varphi_{n} \leq\sup\varphi_{n} < 1$.

\begin{prop}
\emph{\cite[Proposition 21]{LT2}}\label{Proposition0} If $\theta<1$ and
$\inf\varphi_{n}>0,$ then $\left(  \theta_{n}\right)  $ satisfies
\end{prop}

$\left(  \lnot\dag\right)  \qquad$\bigskip\label{Nodagger} $\lim_{m}%
\limsup_{n}\frac{\theta_{m+n}}{\theta_{n}} = 0$ and

$\left(  \ddag\right)  \qquad$\label{doubledagger}There exists
$F:\mathbb{N\rightarrow R}$ with $\lim_{n\rightarrow\infty}F\left(  n\right)
=0$ such that for all $R,t\in\mathbb{N}$ and any arithmetic progression
$\left(  s_{i}\right)  _{i=1}^{R}$ in $\mathbb{N}$,
\[
\max_{1\leq i\leq R}\frac{\theta_{s_{i}+t}}{\theta_{s_{i}}}\leq F\left(
R\right)  \sum_{i=1}^{R}\frac{\theta_{s_{i}+t}}{\theta_{s_{i}}}.
\]

The main tool in our investigation is a construction of certain ``layered
repeated averages" that can be carried out under the assumptions $\left(
\lnot\dag\right)  $ and $\left(  \ddag\right)  .$ The basic units of the
construction are the \emph{repeated averages} due to Argyros, Mercourakis and
Tsarpalias \cite{AMT} which we recall here. An $\mathcal{S}_{0}$-repeated
average is a vector $e_{k}$ for some $k\in\mathbb{N}.$ For any $p\in
\mathbb{N}$, an $\mathcal{S}_{p}$-repeated average is a vector of the form
$\frac{1}{k}\sum_{i=1}^{k}x_{i},$ where $x_{1}<\cdots<x_{k}$ are repeated
$\mathcal{S}_{p-1}$-repeated averages and $k=\min\operatorname*{supp}x_{1}.$
Observe that any $\mathcal{S}_{p}$-repeated average $x$ is a convex
combination of $\left\{  e_{k}:k\in\operatorname*{supp}x\right\}  $ such that
$\left\Vert x\right\Vert _{\infty}\leq\left(  \min\operatorname*{supp}%
x\right)  ^{-1}$ and $\operatorname*{supp}x\in{\mathcal{S}}_{p}$.\newline

\noindent\textbf{Construction of Layered Repeated Averages}

Assume that $\left(  \lnot\dag\right)  $ and $\left(  \ddag\right)  $ hold.
Given $N\in\mathbb{N}$ and $V\in\left[  \mathbb{N}\right]  $, choose sequences
$\left(  p_{k}\right)  _{k=1}^{N}$ and $\left(  L_{k}\right)  _{k=1}^{N}$ in
$\mathbb{N}$, $L_{k}\geq2$, that satisfy the following conditions:

\begin{enumerate}
\item[(A)] $\dfrac{\theta_{p_{M+1}+n}}{\theta_{n}}\leq\frac{\theta_{1}%
}{24N^{2}}\prod_{i=1}^{M}\theta_{L_{i}p_{i}}$ if $0\leq M\leq N-2$ and $n\geq
p_{N}$ (the vacuous product $\prod_{i=1}^{0}\theta_{L_{i}p_{i}}$ is taken to
be $1$),

\item[(B)] $p_{M+1}>\sum_{i=1}^{M}L_{i}p_{i}$ if $0 < M \leq N-2$,

\item[(C)] $F\left(  L_{M+1}\right)  \leq\frac{\theta_{1}}{144N^{2}}%
\prod_{i=1}^{M}\theta_{L_{i}p_{i}}$ if $0<M\leq N-2.$
\end{enumerate}

If $k\in\mathbb{N}$ and $1\leq M\leq N,$ define $r_{M}\left(  k\right)  $ to
be the integer in $\left\{  1,2,...,L_{M}\right\}  $ such that $L_{M}%
|(k-r_{M}\left(  k\right)  ).$ We can construct sequences of vectors
$\mathbf{x}^{0},\dots,\mathbf{x}^{N}$ with the following properties. Let
$\left(  e_{k}\right)  $ be the unit vector basis of $X=T\left[  \left(
{\mathcal{S}}_{n},\theta_{n}\right)  _{n=1}^{\infty}\right]  .$

\begin{enumerate}
\item[$\boldsymbol{(\alpha)}$] $\mathbf{x}^{0}$ is a subsequence of $\left(
e_{k}\right)  _{k\in V}.$

\item[$\boldsymbol{(\beta)}$] Say $\mathbf{x}^{M}=(x_{j}^{M})$ and $m_{j}%
=\min\operatorname*{supp}x_{j}^{M}.$ Then there is a sequence $(I_{k}^{M+1})$
of integer intervals such that $I_{k}^{M+1}<I_{k+1}^{M+1}$, $%
{\displaystyle\bigcup\limits_{k=1}^{\infty}}
I_{k}^{M+1}=\mathbb{N}$ and each vector $x_{k}^{M+1}\in\mathbf{x}^{M+1}$ is of
the form%
\[
x_{k}^{M+1}=\sum_{j\in I_{k}^{M+1}}a_{j}x_{j}^{M},
\]
where $\theta_{r_{M+1}(k)p_{M+1}}\sum_{j\in I_{k}^{M+1}}a_{j}e_{m_{j}}$ is an
$\mathcal{S}_{r_{M+1}(k)p_{M+1}}$-repeated average. Moreover, the sequence
$(a_{j})_{j=1}^{\infty}$ is decreasing.
\end{enumerate}

Each $x_{k}^{M+1}$ is made up of components of diverse complexities. We
analyze it by decomposing it into components of ``pure forms" in the following
manner. We adhere to the notation in $\boldsymbol{(\beta)}$.\newline

\noindent\textbf{``Pure Forms"}
{Given }$1\leq r_{i}\leq L_{i}${, }$1\leq M\leq N-1,${ write }%
\[
x_{k}^{M+1}\left(  r_{M}\right)  =\sum\limits_{\substack{j\in I_{k}%
^{M+1}\\r_{M}\left(  j\right)  =r_{M}}}a_{j}x_{j}^{M}.
\]
{For }$1\leq s<M,${ define}$\ \ ${\ }%
\[
x_{k}^{M+1}\left(  r_{s},...,r_{M}\right)  =%
{\displaystyle\sum\limits_{\substack{j\in I_{k}^{M+1}\\r_{M}\left(  j\right)
\ =r_{M}}}}
a_{j}x_{j}^{M}\left(  r_{s},...,r_{M-1}\right)  .
\]
{If }$1\leq s\leq M,$ {it is clear that }$x_{k}^{M+1}=\sum x_{k}^{M+1}\left(
r_{s},...,r_{M}\right)  ,${ where the sum is taken over all possible values of
}$r_{s},...,r_{M}.$\newline

Given $r_{1},...r_{N},$ write $p\left(  r_{1},...,r_{j}\right)  =\sum
_{i=1}^{j}r_{i}p_{i},$ $1\leq j\leq N.$ Set
\[
\Phi_{k}^{N}=\frac{\theta_{1}}{2}\sum_{r_{1},...,r_{N-1}}\theta_{p\left(
r_{1},...,r_{N}\left(  k\right)  \right)  }\theta_{r_{N}\left(  k\right)
p_{N}}^{-1}\prod_{i=1}^{N-1}\theta_{r_{i}p_{i}}^{-1}L_{i}^{-1}.
\]
If $p\geq N$, define
\[
\Theta_{p}=\Theta_{p}\left(  N\right)  =\max\bigl\{\prod_{i=1}^{N}\theta
_{\ell_{i}}:\ell_{i}\in{\mathbb{N}},\sum_{i=1}^{N}\ell_{i}=p\bigr\}.
\]
The following estimates are crucial for subsequent computations. From here on,
we fix a $k$ satisfying%

\begin{equation}
\label{eq 1}k\geq42N^{2}\prod_{i=1}^{N}L_{i}\theta_{L_{i}p_{i}}^{-1}.
\end{equation}

\begin{prop}
\emph{\cite[Theorem 20; see also the remark following the proof of the
theorem]{LT2}}\label{UppBddx}
\[
\left\Vert x_{k}^{N}\right\Vert \leq\left(  \frac{2}{N}+4\theta_{1}^{-1}%
\sup_{r_{1},\dots,r_{N-1}}\frac{\Theta_{p\left(  r_{1},\dots,r_{N}\left(
k\right)  \right)  }}{\theta_{p\left(  r_{1},\dots,r_{N}\left(  k\right)
\right)  }}\right)  \Phi_{k}^{N}.
\]

\end{prop}

\begin{prop}
\emph{\cite[Corollary 9]{LT2}}\label{ell1x}
\[
\|x_{k}^{N}\left(  r_{1},...,r_{N-1}\right)  \|_{\ell^{1}}\geq\frac{1}%
{2}\theta_{r_{N}\left(  k\right)  p_{N}}^{-1}\prod_{i=1}^{N-1}\theta
_{r_{i}p_{i}}^{-1}L_{i}^{-1}.
\]

\end{prop}

For all $m\in\mathbb{N}$, $z\in c_{00},$ define%
\[
\left\Vert z\right\Vert _{m}=\theta_{m}\sup\left\{  \sum\left\Vert E_{\ell
}z\right\Vert :\left(  E_{\ell}\right)  \text{ is }\mathcal{S}_{m}%
\text{-admissible}\right\}  .
\]

\begin{prop}
\label{Proposition1}Suppose that $x=x_{k}^{N}=\sum_{i=1}^{\ell}b_{i}e_{m_{i}}%
$, $\left(  z_{i}\right)  $ is a normalized block basis of $(e_{k})$ with
$\min\operatorname*{supp}z_{i}=m_{i},$ $q=\sum_{j=1}^{N}L_{j}p_{j},$ and there
exists $K<\infty$ such that $\left\Vert z_{i}\right\Vert _{s}\geq\frac{1}{K}$
for all $1\leq s\leq q,$ $1\leq i\leq\ell$. Let $z=\sum_{i=1}^{\ell}b_{i}%
z_{i}$. Then%
\[
\left\Vert x\right\Vert \leq\left(  \frac{2}{N}+4\theta_{1}^{-1}\sup
_{r_{1},\dots,r_{N-1}}\frac{\Theta_{p\left(  r_{1},\dots,r_{N}\left(
k\right)  \right)  }}{\theta_{p\left(  r_{1},\dots,r_{N}\left(  k\right)
\right)  }}\right)  K\theta_{1}\left\Vert z\right\Vert .
\]

\end{prop}

\begin{proof}
According to Proposition \ref{UppBddx}, it suffices to show that $\left\vert
\left\vert z\right\vert \right\vert \geq\left(  \theta_{1}K\right)  ^{-1}%
\Phi_{k}^{N}.$ For each $1\leq i\leq\ell$, let $\left(  r_{1},...,r_{N-1}%
\right)  $ be the unique $\left(  N-1\right)  $-tuple such that $m_{i}%
\in\operatorname*{supp}x_{k}^{N}\left(  r_{1},...,r_{N-1}\right)  .$ Since
$\left\Vert z_{i}\right\Vert _{t}\geq\frac{1}{K}$ for $t=q-p\left(
r_{1},...,r_{N-1}\right)  ,$ there exists an $\mathcal{S}_{t}$ admissible
family $\mathcal{G}_{i}$ such that $G\subseteq\operatorname*{supp}z_{i}$ for
all $G\in\mathcal{G}_{i}$ and
\begin{equation}
\left\vert \left\vert z_{i}\right\vert \right\vert _{t}=\theta_{t}\sum
_{G\in\mathcal{G}_{i}}\left\vert \left\vert Gz_{i}\right\vert \right\vert
\geq\frac{1}{K}. \label{fdf}%
\end{equation}
We estimate the norm of $z$ by means of a particular tree $\mathcal{T}$.
$\ $If $0\leq n\leq N$ and $\operatorname*{supp}x_{j}^{N-n}\subseteq
\operatorname*{supp}x_{k}^{N},$ let
\[
E_{j}^{n}=\cup\left\{  \operatorname*{supp}z_{i}:m_{i}\in\operatorname*{supp}%
x_{j}^{N-n}\right\}
\]
and
\[
\mathcal{E}^{n}=\left\{  E_{j}^{n}:\operatorname*{supp}x_{j}^{N-n}%
\subseteq\operatorname*{supp}x_{k}^{N}\right\}  .
\]

\noindent By $\boldsymbol{(\beta)}$ in the construction of $\mathbf{x}^{N-n},$
$E_{s}^{n}$ is an $\mathcal{S}_{r_{N-n}\left(  s\right)  p_{N-n}}$-admissible
union of the sets $\left\{  E_{j}^{n+1}:\operatorname*{supp}x_{j}%
^{N-n-1}\subseteq\operatorname*{supp}x_{k}^{N}\right\}  $ . Hence $%
{\textstyle\bigcup_{n=1}^{N}}
\mathcal{E}^{n}$ is an admissible tree so that
\begin{equation}
\operatorname{ord}\left(  E_{j}^{n+1}\right)  =\operatorname{ord}\left(
E_{s}^{n}\right)  +r_{N-n}\left(  s\right)  p_{N-n}\text{ if }E_{j}%
^{n+1}\subseteq E_{s}^{n}. \label{E1}%
\end{equation}
Note that $\operatorname*{supp}x_{j}^{0}$ is a singleton $\left\{
m_{i}\right\}  $ for some $i$ and hence $E_{j}^{N}=\operatorname*{supp}z_{i}.$
It follows from (\ref{E1}) that $\operatorname{ord}\left(  E_{j}^{N}\right)
=p\left(  r_{1},...,r_{N-1}\right)  +r_{N}\left(  k\right)  p_{N},$ where
$\left(  r_{1},...,r_{N-1}\right)  $ is the unique $\left(  N-1\right)
$-tuple determined by $m_{i}.$ Set $\mathcal{E}^{N+1}=\cup_{i=1}^{\ell
}\mathcal{G}_{i}.$ Since $\mathcal{G}_{i}$ is an $\mathcal{S}_{q-p\left(
r_{1},...,r_{N-1}\right)  }$-admissible family with $\cup_{G\in\mathcal{G}%
_{i}}G\subseteq\operatorname*{supp}z_{i}=E_{j}^{N}\in\mathcal{E}^{N},$
$\mathcal{T=\cup}_{n=0}^{N+1}\mathcal{E}^{n}$ is an admissible tree such that
$\operatorname{ord}\left(  G\right)  =q+r_{N}\left(  k\right)  p_{N}$ for each
of the leaves $G$ of $\mathcal{T}.$ By Lemma \ref{Tree},
$\
{\textstyle\bigcup_{i=1}^{\ell}}
\mathcal{G}_{i}$ is $\mathcal{S}_{r_{N}\left(  k\right)  p_{N}+q}%
$-admissible.
Therefore,
\begin{align*}
\left\vert \left\vert z\right\vert \right\vert  &  \geq\theta_{q+r_{N}\left(
k\right)  p_{N}}\sum_{i=1}^{\ell}b_{i}\sum_{G\in\mathcal{G}_{i}}\left\vert
\left\vert Gz_{i}\right\vert \right\vert \\
&  \geq\theta_{q+r_{N}\left(  k\right)  p_{N}}\sum_{i=1}^{\ell}b_{i}\left(
K\theta_{q-p\left(  r_{1},...,r_{N-1}\right)  }\right)  ^{-1}\,\,\,\,\text{by
\eqref{fdf}}\\
&  =\theta_{q+r_{N}\left(  k\right)  p_{N}}\sum_{r_{1},...,r_{N-1}}\left(
K\theta_{q-p\left(  r_{1},...,r_{N-1}\right)  }\right)  ^{-1}\left\vert
\left\vert x\left(  r_{1},...,r_{N-1}\right)  \right\vert \right\vert
_{\ell^{1}}.
\end{align*}
By the regularity of $\left(  \theta_{n}\right)  $, $\theta_{p\left(
r_{1},...,r_{N-1},r_{N}\left(  k\right)  \right)  }\theta_{q-p\left(
r_{1},...,r_{N-1}\right)  }\leq\theta_{q+r_{N}\left(  k\right)  p_{N}}.$
Applying Proposition \ref{ell1x} to the above gives%
\begin{align*}
\left\vert \left\vert z\right\vert \right\vert  &  \geq\frac{\theta
_{q+r_{N}\left(  k\right)  p_{N}}}{K}\sum_{r_{1},...,r_{N-1}}\frac
{\theta_{p\left(  r_{1},...,,r_{N}\left(  k\right)  \right)  }}{\theta
_{q+r_{N}\left(  k\right)  p_{N}}}\left(  \frac{1}{2}\theta_{r_{N}\left(
k\right)  p_{N}}^{-1}\prod_{i=1}^{N-1}\theta_{r_{i}p_{i}}^{-1}L_{i}%
^{-1}\right) \\
&  =\frac{1}{K}\left(  \frac{1}{2}\sum_{r_{1},...,r_{N-1}}\theta_{p\left(
r_{1},...,r_{N}\left(  k\right)  \right)  }\theta_{r_{N}\left(  k\right)
p_{N}}^{-1}\prod_{i=1}^{N-1}\theta_{r_{i}p_{i}}^{-1}L_{i}^{-1}\right) \\
&  =\left(  \theta_{1}K\right)  ^{-1}\Phi_{k}^{N}.
\end{align*}

\end{proof}

We need a few preparatory results in order to exploit the estimate established
in Proposition \ref{Proposition1}.

\begin{lem}
\label{L1}If $\left(  x_{k}\right)  \prec\left(  e_{k}\right)  ,$
$\varepsilon>0$ and $p\in\mathbb{N}$, then there exists $y\in
\operatorname*{span}\left(  x_{k}\right)  ,$ $\left\vert \left\vert
y\right\vert \right\vert =1$ such that $\left\vert \left\vert y\right\vert
\right\vert _{\mathcal{S}_{p}}<\varepsilon.$
\end{lem}

\begin{proof}
Assume to the contrary. There exist $\varepsilon>0,$ $p\in\mathbb{N}$ such
that for all $y\in\operatorname*{span}\left\{  \left(  x_{k}\right)  \right\}
,$ $\left\vert \left\vert y\right\vert \right\vert _{\mathcal{S}_{p}}%
\geq\varepsilon\left\vert \left\vert y\right\vert \right\vert .$ On the other
hand, $\left\vert \left\vert y\right\vert \right\vert \geq\theta_{p}\left\vert
\left\vert y\right\vert \right\vert _{\mathcal{S}_{p}}.$ Hence $\left\vert
\left\vert \cdot\right\vert \right\vert $ and $\left\vert \left\vert
\cdot\right\vert \right\vert _{\mathcal{S}_{p}}$ are equivalent on
$\operatorname*{span}\left\{  \left(  x_{k}\right)  \right\}  .$ However, the
Schreier space $\mathcal{S}_{p}$ is $c_{0}$-saturated. It follows that
$\left[  \left(  x_{k}\right)  \right]  $ and thus $X$ contains a copy of
$c_{0},$ contradicting the reflexivity of $X$.
\end{proof}

\begin{lem}
\label{L3}If $\left(  z_{k}\right)  \prec\left(  y_{k}\right)  \prec\left(
e_{k}\right)  $, and $\left\vert \left\vert y_{k}\right\vert \right\vert
_{\mathcal{S}_{k-1}}\leq\frac{1}{2^{k+2}}$ for all $k,$ then
\[
\left\vert \left\vert z_{k}\right\vert \right\vert _{\mathcal{S}_{k-1}}%
\leq\frac{1}{2^{k+1}}\text{ for all }k.
\]

\end{lem}

\begin{proof}
Write $z_{k}=\sum_{j\in J_{k}}a_{j}y_{j}$. Note that $\left\vert
a_{j}\right\vert \leq\left\vert \left\vert z_{k}\right\vert \right\vert =1$
for all $j\in J_{k}.$ Therefore,
\begin{align*}
\left\vert \left\vert z_{k}\right\vert \right\vert _{\mathcal{S}_{k-1}}  &
\leq\sum_{j\in J_{k}}\left\vert \left\vert y_{j}\right\vert \right\vert
_{\mathcal{S}_{k-1}}\\
&  \leq\sum_{j\in J_{k}}\left\vert \left\vert y_{j}\right\vert \right\vert
_{\mathcal{S}_{j-1}}\text{ since }k\leq\min J_{k}\leq j\\
&  \leq\sum_{j\in J_{k}}\frac{1}{2^{j+2}}\leq\frac{1}{2^{k+1}}.
\end{align*}

\end{proof}

\begin{lem}
\label{L2}Assume that $\theta<1$ and $\inf_{n}\varphi_{n}>0.$ $\ $If $\left(
z_{k}\right)  \prec\left(  e_{k}\right)  $ and
\[
\left\vert \left\vert z_{k}\right\vert \right\vert _{\mathcal{S}_{k-1}}%
\leq\frac{1}{2^{k+1}}\text{ for all }k,
\]
then there is a constant $K<\infty$ such that for all $z\in
\operatorname*{span}\left(  z_{k}\right)  _{k=n}^{\infty},$ we have
$\left\vert \left\vert z\right\vert \right\vert _{m}\geq\frac{1}{2K}\left\vert
\left\vert z\right\vert \right\vert $ for all $1\leq m\leq n.$
\end{lem}

\begin{proof}
First observe that
\begin{equation}
\label{fdf2}\frac{\theta_{m+n}}{\theta_{m}\theta_{n}}=\frac{\varphi_{m+n}%
}{\varphi_{m}\varphi_{n}}\leq\frac{1}{\left(  \inf\varphi_{n}\right)  ^{2}%
}\text{ for all }m,n.
\end{equation}
Let $K=\frac{1}{\left(  \inf\varphi_{n}\right)  ^{2}}.$ Suppose that
$z\in\operatorname*{span}\left(  z_{k}\right)  _{k=n}^{\infty},$ $\left\vert
\left\vert z\right\vert \right\vert =1$ and $1\leq m\leq n.$ Choose an
admissible tree $\mathcal{T}$ of $z$ so that
\begin{align*}
1  &  =\left\vert \left\vert z\right\vert \right\vert =\mathcal{T}z=\sum
_{E\in\mathcal{L}\left(  \mathcal{T}\right)  }t\left(  E\right)  \left\vert
\left\vert Ez\right\vert \right\vert _{c_{0}}\\
&  =\sum_{\substack{E\in\mathcal{L}\left(  \mathcal{T}\right)
\\\operatorname{ord}\left(  E\right)  \leq m}}t\left(  E\right)  \left\vert
\left\vert Ez\right\vert \right\vert _{c_{0}}+\sum_{\substack{E\in
\mathcal{L}\left(  \mathcal{T}\right)  \\\operatorname{ord}\left(  E\right)
>m}}t\left(  E\right)  \left\vert \left\vert Ez\right\vert \right\vert
_{c_{0}}.
\end{align*}
Write $z=\sum_{k=n}^{\infty}a_{k}z_{k}.$ Then $\left\vert a_{k}\right\vert
\leq1$ as $\left\vert \left\vert z\right\vert \right\vert =1.$ Note that
according to Lemma \ref{Tree},
the collection $\left\{  E\in\mathcal{L}\left(  \mathcal{T}\right)
:\operatorname{ord}\left(  E\right)  \leq m\right\}  $ of leaves is
$\mathcal{S}_{m}$-admissible. Therefore,%
\begin{align*}
\sum_{\substack{E\in\mathcal{L}\left(  \mathcal{T}\right)
\\\operatorname{ord}\left(  E\right)  \leq m}}t\left(  E\right)  \left\vert
\left\vert Ez\right\vert \right\vert _{c_{0}}  &  \leq\sum_{\substack{E\in
\mathcal{L}\left(  \mathcal{T}\right)  \\\operatorname{ord}\left(  E\right)
\leq m}}\left\vert \left\vert Ez\right\vert \right\vert _{c_{0}}\\
&  \leq\left\vert \left\vert z\right\vert \right\vert _{\mathcal{S}_{m}}\\
&  \leq\sum_{k=n}^{\infty}\left\vert \left\vert z_{k}\right\vert \right\vert
_{\mathcal{S}_{m}}\\
&  \leq\sum_{k=n}^{\infty}\left\vert \left\vert z_{k}\right\vert \right\vert
_{\mathcal{S}_{k-1}}\leq\sum_{k=n}^{\infty}\frac{1}{2^{k+1}}\leq\frac{1}{2}.
\end{align*}
Thus
\[
\sum_{\substack{E\in\mathcal{L}\left(  \mathcal{T}\right)
\\\operatorname{ord}\left(  E\right)  >m}}t\left(  E\right)  \left\vert
\left\vert Ez\right\vert \right\vert _{c_{0}}\geq\frac{1}{2}.
\]
Let $\mathcal{E}$ be the collection of all nodes $E$ in $\mathcal{T}$ that are
minimal subject to the condition $\operatorname{ord}\left(  E\right)  >m$.
Also, let $\mathcal{D}$ be the set of all immediate predecessors of nodes in
$\mathcal{E}$. If $D\in\mathcal{D}$, let $\mathcal{E}\left(  D\right)  $ be
the collection of its immediate successors. For each $E\in\mathcal{E}\left(
D\right)  ,$ $\operatorname{ord}\left(  D\right)  \leq m<\operatorname{ord}%
\left(  E\right)  .$ Therefore there exists an $\mathcal{S}%
_{m-\operatorname{ord}\left(  D\right)  }$-admissible collection
$\mathcal{G}_{D}$ of subsets of $D$ such that $\mathcal{E}\left(  D\right)
=\cup\left\{  E\in\mathcal{E}\left(  D\right)  :E\subseteq G\text{ for some
}G\in\mathcal{G}_{D}\right\}  $ and $\left\{  E\in\mathcal{E}\left(  D\right)
:E\subseteq G\right\}  $ is $\mathcal{S}_{\operatorname{ord}\left(  E\right)
-m}$-admissible for each $G\in\mathcal{G}_{D}.$ Now $\mathcal{G}=\cup
_{D\in\mathcal{D}}\mathcal{G}_{D}$ is $\mathcal{S}_{m}$-admissible and
$\theta_{\operatorname{ord}\left(  E\right)  }\geq t\left(  E\right)  $ by the
regularity of $\left(  \theta_{n}\right)  $. Hence
\begin{align*}
\left\vert \left\vert z\right\vert \right\vert _{m}  &  \geq\theta_{m}%
\sum_{G\in\mathcal{G}}\left\vert \left\vert Gz\right\vert \right\vert \\
&  \geq\theta_{m}\sum_{G\in\mathcal{G}}\theta_{\operatorname{ord}\left(
E\right)  -m}\sum_{\substack{E\in\mathcal{E}\\E\subseteq G}}\left\vert
\left\vert Ez\right\vert \right\vert \\
&  \geq\sum_{G\in\mathcal{G}}\frac{\theta_{\operatorname{ord}\left(  E\right)
}}{K}\sum_{\substack{E\in\mathcal{E}\\E\subseteq G}}\left\vert \left\vert
Ez\right\vert \right\vert \text{ by \eqref{fdf2} and the definition of }K\\
&  \geq\frac{1}{K}\sum_{E\in\mathcal{E}}t\left(  E\right)  \left\vert
\left\vert Ez\right\vert \right\vert \geq\frac{1}{K}\sum_{\substack{E\in
\mathcal{L}\left(  \mathcal{T}\right)  \\\operatorname{ord}\left(  E\right)
>m}}t\left(  E\right)  \left\vert \left\vert Ez\right\vert \right\vert
\geq\frac{1}{2K}.
\end{align*}

\end{proof}

We shall show that, for appropriate $(\theta_{n})$, if $\left(  z_{k}\right)
\prec\left(  e_{k}\right)  $ satisfies the conclusion of Lemma \ref{L3}, then
it is not equivalent to a subsequence of $\left(  e_{k}\right)  .$

\begin{lem}
\label{L4}If $0<\inf_{n}\varphi_{n}\leq\sup_{n}\varphi_{n}<1,$ then $\lim
_{N}\sup_{p\geq N}\frac{\Theta_{p}\left(  N\right)  }{\theta_{p}}=0.$
\end{lem}

\begin{proof}
Let $\varepsilon>0.$ Choose $N$ such that $\frac{d^{N}}{c}<\varepsilon,$ where
$0<c=\inf_{n}\varphi_{n}\leq\sup_{n}\varphi_{n}=d<1.$ Let $p\in\mathbb{N}$. If
$\left(  \ell_{i}\right)  _{i=1}^{N}$ is a sequence of positive integers such
that $\sum_{i=1}^{N}\ell_{i}=p,$ then
\[
\prod_{i=1}^{N}\theta_{\ell_{i}}=\theta^{p}\prod_{i=1}^{N}\varphi_{\ell_{i}%
}\leq\theta^{p}d^{N}%
\]
and
\[
\theta_{p}=\varphi_{p}\theta^{p}\geq c\theta^{p}.
\]
Thus
\[
\sup_{p\geq N}\frac{\Theta_{p}\left(  N\right)  }{\theta_{p}}\leq\frac{d^{N}%
}{c}<\varepsilon.
\]

\end{proof}

\begin{prop}
\label{Proposition12}If $\left(  z_{k}\right)  $ is a normalized block basis
that is equivalent to a subsequence of $\left(  e_{k}\right)  ,$ then there is
a subsequence $\left(  z_{k_{j}}\right)  $ of $\left(  z_{k}\right)  $ such
that $\left(  z_{k_{j}}\right)  $ is equivalent to $\left(  e_{m_{j}}\right)
,$ where $m_{j}=\min\operatorname*{supp}z_{k_{j}}.$
\end{prop}

\begin{proof}
It is well known that two subsequences $\left(  e_{n_{i}}\right)  $ and
$\left(  e_{\ell_{i}}\right)  $ of $\left(  e_{k}\right)  $ are equivalent
whenever $\max\left\{  n_{i},\ell_{i}\right\}  <\min\left\{  n_{i+1}%
,\ell_{i+1}\right\}  $ for all $i.$ If $\left(  z_{k}\right)  $ is equivalent
to a subsequence of $\left(  e_{k}\right)  ,$ then there is a subsequence
$\left(  z_{k_{j}}\right)  $ of $\left(  z_{k}\right)  $ that is equivalent to
a subsequence $\left(  e_{n_{j}}\right)  $ of $\left(  e_{k}\right)  $ with
\[
\max\{\min\operatorname*{supp}z_{k_{j}},n_{j}\}<\min\left\{  \min
\operatorname*{supp}z_{k_{j+1}},n_{j+1}\right\}  \text{ for all }j.
\]
Thus $\max\left\{  n_{j},m_{j}\right\}  <\min\left\{  n_{j+1},m_{j+1}\right\}
,$ and hence $\left(  e_{n_{j}}\right)  $ is equivalent to $\left(  e_{m_{j}%
}\right)  .$ Consequently, $\left(  z_{k_{j}}\right)  $ is equivalent to
$\left(  e_{m_{j}}\right)  .$
\end{proof}

We are now ready to prove the main result of the section.

\begin{thm}
\label{Theorem1}If $\ ~0<\inf_{n}\varphi_{n}\leq\sup_{n}\varphi_{n}<1,$ then
$X$ is strongly non-subsequentially minimal.
\end{thm}

\begin{proof}
Let $\left(  x_{k}\right)  $ be a normalized block basis of $\left(
e_{k}\right)  .$ By Lemma \ref{L1}, there exists $\left(  y_{k}\right)
\prec\left(  x_{k}\right)  $ such that $\left\Vert y_{k}\right\Vert
_{\mathcal{S}_{k-1}}\leq\frac{1}{2^{k+2}}$ for all $k.$ Suppose that there
exists $\left(  z_{k}\right)  \prec\left(  y_{k}\right)  $ that is equivalent
to a subsequence of $\left(  e_{k}\right)  $. Applying Proposition
\ref{Proposition12}, we may assume that $\left(  z_{k}\right)  $ is equivalent
to $\left(  e_{m_{k}}\right)  ,$ where $m_{k}=\min\operatorname*{supp}z_{k}.$
Pick $\varepsilon>0$ so that
\[
\varepsilon\left\Vert
{\displaystyle\sum}
b_{k}z_{k}\right\Vert \leq\left\Vert
{\displaystyle\sum}
b_{k}e_{m_{k}}\right\Vert \text{ for all }\left(  b_{k}\right)  \in c_{00}.
\]
By a combination of Lemmas \ref{L3} and \ref{L2} there is a constant
$K<\infty$ such that $\left\vert \left\vert z\right\vert \right\vert _{s}%
\geq\frac{1}{2K}\left\vert \left\vert z\right\vert \right\vert ,$ for all
$z\in\operatorname*{span}\left(  z_{k}\right)  _{k=n}^{\infty},$ $1\leq s\leq
n.$\bigskip\ Use Lemma \ref{L4} to choose $N$ such that $\frac{2}{N}%
+4\theta_{1}^{-1}\sup_{p}\frac{\Theta_{p}\left(  N\right)  }{\theta_{p}}%
<\frac{\varepsilon}{2K\theta_{1}}$ if $p\geq N.$ With the chosen $N$ and
$V=\left(  m_{i}\right)  _{i=q}^{\infty}$ construct the layered repeated
average vector $x=x_{k}^{N}=\sum_{i=q}^{\ell}b_{i}e_{m_{i}}$ with $k$
satisfying the inequality (\ref{eq 1}). Let $z=\sum_{i=q}^{\ell}b_{i}z_{i}.$
(Recall that $q=\sum_{j=1}^{N}L_{j}p_{j},$ where $\left(  p_{j}\right)
_{j=1}^{N}$ and $\left(  L_{j}\right)  _{j=1}^{N}$ are chosen to satisfy
conditions \textbf{(A)}, \textbf{(B)}, and \textbf{(C)} once $N$ is
determined.) According to Proposition \ref{Proposition1},
\begin{align*}
\left\Vert x\right\Vert  &  \leq\left(  \frac{2}{N}+4\theta_{1}^{-1}%
\sup_{r_{1},\dots,r_{N-1}}\frac{\Theta_{p\left(  r_{1},\dots,r_{N}\left(
k\right)  \right)  }}{\theta_{p\left(  r_{1},\dots,r_{N}\left(  k\right)
\right)  }}\right)  2K\theta_{1} \left\Vert z\right\Vert \\
&  <\varepsilon\left\Vert z\right\Vert ,
\end{align*}
contrary to the choice of $\varepsilon.$
\end{proof}

The following example shows that the condition $\sup_{n}\varphi_{n}<1$ is not
necessary for the conclusion of the theorem to hold.

\begin{example}
\label{Ex2}If $\theta<1,$ there exists a regular sequence $\left(  \theta
_{n}\right)  $ with $\sup_{n}\theta_{n}^{1/n}=\theta$ and $\lim_{n}\varphi
_{n}=1$ such that $X$ is strongly non-subsequentially minimal.
\end{example}

\begin{proof}
Suppose that $0<\theta<1.$ In \cite[Example 23]{LT2}, a regular sequence
$\left(  \theta_{n}\right)  $ is constructed so that $\sup_{n}\theta_{n}%
^{1/n}=\theta,$ $\lim_{n}\varphi_{n}=1$ and for all $N\in\mathbb{N},$ there
are sequences $\left(  p_{k}\right)  _{k=1}^{N}$ and $\left(  L_{k}\right)
_{k=1}^{N}$ satisfying conditions \textbf{(A)}, \textbf{(B)}, and \textbf{(C)}
and
\begin{equation}
\lim_{N\rightarrow\infty}\sup_{r_{1},\dots,r_{N-1}}\frac{\Theta_{p\left(
r_{1},\dots,r_{N}\left(  k\right)  \right)  }}{\theta_{p\left(  r_{1}%
,\dots,r_{N}\left(  k\right)  \right)  }}=0. \label{E2}%
\end{equation}
Following the arguments in Theorem \ref{Theorem1} with Lemma \ref{L4} replaced
by (\ref{E2}) shows that $X$ is strongly non-subsequentially minimal.
\end{proof}

In view of Proposition \ref{qm}, any subsequentially minimal partly modified
mixed Tsirelson space is quasi-minimal. However, the existence of strongly
non-subsequentially minimal mixed Tsirelson spaces prompts the following
question.\newline

\noindent\textbf{Question.} Does every (partly modified) mixed Tsirelson space
$T[({\mathcal{S}}_{n}, \theta_{n})^{\infty}_{n=1}]$ (or $T[({\mathcal{S}}_{n},
\sigma_{n}, \theta_{n})^{\infty}_{n=1}]$) contain a quasi-minimal subspace?

\end{document}